\documentclass{article}

\usepackage{amsmath}
\usepackage{amsfonts}
\usepackage{amssymb}
\usepackage{color}  
\usepackage{amsthm}
\usepackage{graphicx}
\usepackage{caption}
\usepackage{comment}
\usepackage{subcaption}
\usepackage{algpseudocode}
\usepackage{algorithm}
\usepackage[utf8]{inputenc}
\usepackage{tabularx}
\usepackage{multirow}
\usepackage{hyperref}
\usepackage[english]{babel}
\usepackage{setspace}

\usepackage{natbib}
\bibpunct[, ]{(}{)}{,}{a}{}{,}%

\title{Local Area Routes for Vehicle Routing Problems}

\author{Udayan Mandal\textsuperscript{ \rm1}, Amelia Regan\textsuperscript{ \rm 1}, Julian Yarkony \textsuperscript{\rm 1, \rm 2} \\
\textsuperscript{\rm 1}University of California, Irvine, CA\\
\textsuperscript{\rm 2} Laminaar Optimization Research Group, La Jolla, CA
}\date{July 2022}
\begin{document}
\maketitle
\begin{abstract}
We consider an approach for improving the efficiency of column generation (CG) methods for solving vehicle routing problems.  We introduce Local Area (LA) route relaxations, an alternative/complement to the commonly used ng-route relaxations and Decremental State Space Relaxations (DSSR) inside of CG formulations. LA routes are a subset of ng-routes and a super-set of elementary routes. Normally, the pricing stage of CG must produce elementary routes, which are  routes without repeated customers, using processes which can be computationally expensive. Non-elementary routes visit at least one customer more than once, creating a cycle. LA routes relax the constraint of being an elementary route in such a manner as to permit efficient pricing. 
LA routes are best understood in terms of ng-route relaxations. Ng-routes are routes which are permitted to have non-localized cycles in space; this means that at least one intermediate customer (called a breaker) in the cycle must consider the starting customer in the cycle to be spatially far away. LA routes are described using a set of special indexes corresponding to customers on the route ordered from the start to the end of the route.  LA route relaxations further restrict the set of permitted cycles beyond that of ng-routes by additionally enforcing that the breaker must be a located at a special index where the set of special indexes is defined recursively as follows. The first special index in the route is at index 1 meaning that it is associated with the first customer in the route. The k'th special index corresponds to the first customer after the k-1'th special index, that is not considered to be a neighbor of (considered spatially far from) the customer located at the k-1'th special index.
We demonstrate that LA route relaxations can significantly improve the computational speed of pricing when compared to the standard DSSR.
\end{abstract}

\section{Introduction}
In this document we introduce a new tool called Local Area routes (LA routes), which serves as a component in exact column generation (CG) solutions \citep{barnprice,Desrochers1992,cuttingstock} to vehicle routing problems.  For purposes of exposition for LA routes, we define the Capacitated Vehicle Routing Problem (CVRP) below, which is a specific vehicle routing problem, though our approach is not limited to such problems.  CVRP defines a problem consisting of a starting and ending depot, a set of customers with integer demands, and a set of homogeneous vehicles with integer capacity. The customers and the starting/ending depot are positioned at various locations in two-dimensional space. Each vehicle starts and ends at the depot, and does not service more demand than its capacity. We select a set of routes so as to minimize total travel distance (sum of the travel distances of the individual routes) subject to the constraint that all customers are serviced.

Exact solutions to vehicle routing problems are traditionally formulated as weighted set cover problems where each customer must be covered and routes describe feasible sets with the cost of a set being the cost of the route.  This relaxation is a much tighter relaxation than that of compact forms\citep{Desrochers1992,costa2019}.  Since the number of such routes grows exponentially in the number of customers, CG methods are applied so that the set of routes need not be explicitly enumerated.  CG methods imitate the revised simplex approach and generate primal variables during pricing (called columns) on demand, which is tackled as a combinatorial optimization problem. 

In such problems, solving the CG pricing problem to generate routes for the CG restricted master problem (RMP) is especially problematic since the routes produced must be elementary (meaning that no customer is visited more than once in the route). We seek to ease solving the pricing problem by producing a class of routes known as Local Area (LA) routes that are easy to price over, and can be used to efficiently generate elementary routes when used inside Decremental State Space Relaxation (DSSR) \citep{righini2008new,righini2009decremental}.


%
An LA route is a route in which elementarity is relaxed, but not by much due to other constraints. The set of LA routes does not include routes with cycles localized in space; where a cycle is a section of a route consisting of the same customer at the start and end of the section. Localized cycles in space are cycles consisting of customers which are all spatially close to one another. 
LA routes are a subset of the popular ng-routes \citep{baldacci2011new} (and a superset of elementary routes). Ng-routes are a part of many efficient and modern CG algorithms for vehicle routing problems \citep{costa2019}.  Both LA routes and ng-routes can be used inside DSSR to produce exact solutions to pricing over elementary routes.  

No optimal integer solution of the set cover formulation over non-elementary routes that are otherwise feasible (obey capacity restriction) visits the same customer more than once (either within a route or over multiple routes) since by removing that customer we could decrease the cost of the solution, while preserving feasibility.  However, an optimal fractional solution may use such routes.  Hence solving optimization over a more restrictive subset of routes is beneficial if we cannot use elementary routes for computational reasons. Solving the set cover formulation over a class of relaxed routes can be used inside a branch-and-price formulation to ensure an exact solution to the set cover formulation \citep{barnprice}.  
%
 \subsection{Types of Routes}
We now consider various classes of routes that are defined in the CVRP literature in addition to LA routes. To assist in this discussion we consider the following graph upon which such routes are defined. Consider a directed acyclic graph $G$ with vertex set $V$ and edge set $E$.  Vertices are indexed by $i$ or $j$ and edges are indexed by $ij$. We define the starting depot to be the source, denoted by $-1$ and the ending depot to be the sink, denoted by $-2$. The starting/ending depot are the same place but we treat them separately for convenience.  

Each node except $(-1,-2)$ is associated with a customer and the remaining capacity in the vehicle. A path reaching node $i=(u,d)$ indicates that prior to servicing $u$ that there are $d$ units of capacity remaining in the vehicle. We define $d_u$ to be the demand required by customer $u$. 
We connect $i=(u,d)$ to $j=(v,d-d_u)$ for each pair of customers $u,v$ s.t. $u\neq v,d-d_u\geq d_v$ to form an edge with weight equal to the distance from $u$ to $v$. Traversing this edge indicates that the vehicle leaves $u$ with $d-d_u$ units of capacity remaining and travels immediately to customer $v$.  
We connect $-1$ to $(u,d_0)$, where $d_0$ is the vehicle capacity, to form an edge with weight equal to the distance from the depot to $u$. Traversing this edge indicates that $u$ is the first customer visited on the route.  We connect each node $(u,d)$ (for $d_u\leq d$) to the sink $-2$ to form edges with weights equal to the distance from $u$ to the depot. Traversing any one of these edges indicates that after servicing $u$ the vehicle heads to the ending depot, terminating its route. Each elementary path starting at the source and ending at the sink corresponds to an elementary route. A path is elementary if there is no more than one node corresponding to a given customer on the path. Each path has an associated cost equal to the sum of weights of edges on the path.

We use $P$ to denote the set of paths starting at $-1$ and ending at $-2$.  Each such path in $P$ describes a route that may not be elementary but does not service more demand than available capacity and starts and ends at the depot. Such paths are referred to as resource feasible. We use $V_p$ to denote the set of nodes excluding the source and sink in path $p$. 
We use $u^p_{k}$ to denote the customer associated with the $k$'th node excluding the source and the sink in path $p$.  We now describe various classes of routes.  

\begin{itemize}
\item \textbf{Elementary Routes:  }

We define $\Omega$ to be the set of all elementary routes, which are paths in $P$ that visit no customer in the path more than once. We define $\Omega$ formally using the following equation, with exposition afterwards.

\begin{align}
   \Omega=\{ p \in P; u^p_{k_1}\neq u^p_{k_2} \quad \forall k_1,k_2 \quad s.t. (0<k_1<k_2\leq |V_p|) 
\end{align}

$\Omega$ consists of any path in $P$ that does not visit the same customer at $k_1$ and $k_2$, where $k_1$ and $k_2$ ($k_1$ comes before $k_2$) are indexes of nodes on the path excluding the source and the sink.

\item \textbf{Q-routes:  }

Q-routes are routes which can not have cycles of length $1$. This means that a route visiting customer 1, followed by customer 2, and then visiting customer 1 again is forbidden. However, a route visiting customer 1, then visiting customer 2, then visiting customer 3, and then visiting customer 1 after is not forbidden. Q-routes were introduced by Christofides, Mingozzi and Toth in order to aid CG in solving vehicle routing problems \citep{christofides1981exact}.

We use $\Omega^1$ to denote the set of Q-Routes, which is defined as follows. 
\begin{align}
    \Omega^1=\{ p \in P; u_{k}\neq u_{k+2} \quad \forall k \quad s.t. (0<k, k+2\leq |V_p|) 
\end{align}
Q-routes can be generalized as KQ-routes where KQ-routes enforce that a customer can only be visited again after visiting at least $K$ intermediate customers. Observe that Q-routes correspond to KQ-routes with K=1. We use the term KQ to draw the link to Q-routes, but this terminology is not standard in the literature.  

We use $\Omega^K$ to denote the set of KQ-Routes, which is defined as follows. 
\begin{align}
   \Omega^K=\{ p \in P; u^p_{k_1}\neq u^p_{k_2} \quad \forall k_1,k_2 \quad s.t.  (0<k_1<k_2\leq \min(|V_p|,k_1+K)) 
\end{align}
 
\item \textbf{ng-routes:  }

Ng-routes are highly celebrated and used by many researchers \citep{baldacci2011new}. Each customer is associated with a set of customers which are close in proximity to that customer (also known as neighbors of that customer). This set of neighbors of $u$ is denoted $M_u$, where $u$ represents the customer the set is associated with. Ng-routes ban spatially localized cycles by enforcing that a cycle can only exist starting and ending at $u$ if there is an intermediate customer $v$ for which $u\notin M_v$.   

\item \textbf{LA-routes:}

Local Area routes (LA routes) are the topic of this paper. LA routes are a subset of ng-routes but further restrict cycles. Thus the CVRP set cover LP relaxation over LA routes is no looser than, and in fact potentially tighter than that over ng-routes. LA routes are defined using LA neighborhood set $N_u$ (for each customer $u$) where $N_u$ consists of spatially nearby customers to $u$.  The LA neighborhood sets $N_u$ are computationally easier to consider than $M_u$ and hence can be larger than $M_u$.  LA routes are defined with a set of special indexes associated with each path $p$.  
%
The set of special indexes in a path $p$ is defined recursively from the start of the route, with the first special index being equal to one.  Let $q^p_j$ be the index of the $j$'th special index (meaning $q^p_1=1$).  
The $j$'th special index corresponds to the first customer after $q^p_{j-1}$ that is not considered to be the an LA neighbor of $u^p_{q^p_{j-1}}$. 
We define the set of special indexes by defining the $q^p_j$ terms recursively as follows.  
\begin{align}
    q_{j}^p\leftarrow \min_{\substack{k>q_{j-1}^p\\ u^p_k \notin N_v}} k \quad \forall j>1; \mbox{where $u^p_{q_{j-1}}=v$}\\
    q_1^p=1 \nonumber 
\end{align}
A resource feasible path (excluding elementarity) is an LA-route if for any cycle in that path starting/ending at customer $u$, there is an intermediate special index with associated customer $v$ for which $u \notin M_v$ (note the use of $M_v$ not $N_v$ here).  Note that the difference between an ng-route and an LA route is understood by building on the statement ``Ng-routes ban spatially localized cycles by enforcing that a cycle can only exist starting and ending at $u$ if there is an intermediate customer $v$ for which $u\notin M_v$."  In contrast LA-routes ban spatially localized cycles by enforcing that a cycle can only exist starting and ending at $u$ if there is an intermediate customer $v$ {\color{red} at a special index} for which $u\notin M_v$. 
\end{itemize}

\subsection{Overview of Efficient Pricing over LA Routes  }

LA routes possess mathematical properties that make pricing over them efficient. We now provide an overview of efficient pricing over LA routes, which exploits the following three properties. 
\begin{enumerate}
    \item For each $u \in N, w\in \hat{N}, v\in \hat{N}, \hat{N}\subseteq N_u$ (where $N$ is the set of all customers) we can compute the lowest cost elementary path starting at $w$ ending at $v$ and visiting all customers in $\hat{N}$. The order of customers in this path does not change as a function of the dual variables. The size of the sets $N_u$ is selected to be small enough (10-20) to permit all such terms to be computed in advance of CG (not at each iteration of CG).  This is done as a dynamic program and detailed in Section \ref{Sec_LA_fast}.
    \item Given any $u\in N$, $v \in (N\cup -2)-(N_u\cup u), \hat{N}\subseteq N_u$  we can compute the lowest cost elementary path starting at $u$ ending at $v$ and servicing all customers in $\hat{N}$ using the computed terms in bullet point 1. The order of customers in this path does not change as a function of the dual variables. This is computed once at the beginning of CG, not at each iteration of CG.  
    \item Given any $u \in N, v \in (N\cup -2)-(N_u\cup u), d \geq d_u$ we can at each iteration of CG during pricing, 
    efficiently compute the lowest reduced cost elementary path starting at $u$, ending at $v$, servicing only intermediate customers in $N_u$ and servicing $d$ units of demand (including $u$ not $v$).  The membership and the ordering (given the membership, which is the customers included in the path) of this path does not change as a function of the dual variables (by bullet point 2).
\end{enumerate}  

We now discuss briefly how pricing is done over LA routes when the ng-neighborhood sets $M_u$ are empty, and then we extend this to the case where $M_u$ sets are not empty. When the $M_u$ sets are not empty, we can integrate LA route relaxations with DSSR. Greater exposition is provided in later sections of this document. To assist in this discussion we introduce a directed acyclic graph with nodes identical to the previous graph. As in our previous graph each node except $(-1,-2)$ is associated with a customer and the remaining capacity in the vehicle. As before, a path reaching node $i=(u,d)$ indicates that prior to servicing $u$ there are $d$ units of capacity remaining in the vehicle. 
We connect $i=(u,d_1)$ to $j=(v,d_2)$ for $u\in N,v\in N-(N_u \cup u)$, s.t. $d_1-d_u \geq d_2\geq d_v$ to form an edge with weight equal to the reduced cost of the lowest reduced cost elementary path servicing exactly $d_1-d_2$ units of demand (including $u$, not $v$), starting at $u$ and ending at $v$, and where all intermediate customers (denoted $N_{ij}$) lie in $N_u$.  Traversing this edge indicates that the vehicle leaves $u$ with $d_1-d_u$ units of capacity remaining and travels to customer $v$ servicing $d_1-d_u-d_2$ units of intermediate demand, and visits the intermediate customers $N_{ij}$ in the optimal order (where optimality corresponds to minimizing the total cost). Note that the optimal ordering of customers for any such edge is not a function of $d_1$,$d_2$ but only $d_1-d_2$ and hence does not need to be computed for each $d_1,d_2$ combination. This edge is nonexistent if no such elementary path exists.  In the remainder of this section when no such path exists then the edge is not created, in this and subsequent cases.    

We connect $-1$ to $(u,d_0)$, where $d_0$ is the vehicle capacity, to form an edge with weight equal to the distance from the depot to $u$. Traversing this edge indicates that $u$ is the first customer visited on the route. We connect each node $(u,d)$ (for $d_u\leq d$) to the sink $-2$ to form an edge with weight equal to the reduced cost of the lowest reduced cost elementary path starting at $u$, ending at the depot, and servicing up to $d$ units of demand (including $u$), where all intermediate customers (denoted $N_{i,-2}$) lie in $N_u$. Traversing any one of these edges indicates that after servicing $u$, the vehicle heads to the ending depot, terminating its route and servicing all customers in $N_{i,-2}$ in the optimal order.

The lowest reduced cost elementary route can be represented by a path in this graph, as well as other paths which need not be of lowest reduced cost. However, large numbers of paths containing localized cycles in space can not be represented and the set of feasible paths decreases as the sizes of $N_u$ sets increase. All paths from the source to the sink correspond to LA routes when all $M_u$ sets are empty. In order to permit the use of Dijkstra's algorithm for computing the lowest reduced cost path, which is computationally faster than the Bellman-Ford algorithm, edge weights need to be non-negative.  Thus we offset edge weights by constant $\eta$ weighted by the the change in demand between nodes (with nodes $-1,-2$ having demands $0$ respectively) serviced on the edge so as to ensure weights are non-negative. Thus we add $\eta*(d_1-d_2)$ to the edge weight between $(u,d_1)$ and $(v,d_2)$ where the source and sink are referenced as $(-1,d_0),(-2,0)$ respectively.  Here $-\eta$ is defined to be the minimizer of the reduced cost of the edge divided by the amount of demand serviced (over all edges), thus $-\eta*d_0$ is a lower bound on the reduced cost of the route.  By subtracting $\eta d_0$ from the cost of any path we get the associated reduced cost.  The optimal route is not changed by this alteration to the weights, since all paths have costs offset by a constant value $\eta d_0$.  

DSSR can then be applied to the $M_u$ neighbor sets (not the $N_u$ neighbor sets) in order to enforce that the route generated above is elementary. In such cases, nodes become associated with $(u,M_1,d)$ terms, where $M_1 \subseteq M_u$ and being at node $u,M_1,d$ indicates that the nascent route is at customer $u$ with $d$ units of demand remaining (prior to servicing $u$) and has serviced each customer in $M_1$ at least once. Edge weights are suitably adjusted given edges between nodes.  For example, we connect each node $(u,M_1,d)$ (for $d_u\leq d$) to the sink $-2$ to form edges with weights equal to the reduced cost of the lowest reduced cost elementary path starting at $u$, ending at the depot, and servicing up to $d$ units of demand (including $u$) and where all intermediate customers lie in $N_u-M_1$. We connect each node $(u,M_1,d_1)$ to another node $(v,M_2,d_2)$ with edge weight corresponding to the cost of the lowest reduced cost path starting at $u$, ending at $v$, consuming demand $d_1 - d_2$ (excluding $v$ and including $u$), servicing some customers in the set $N_u - M_1$, servicing no customers in the set $M_v - M_2$, and servicing all customers in the set $M_2 - M_1$. More about these connections are described later in the document in Section \ref{Sec_la_ILP}.

%

To assist DSSR we use $A*$ \citep{dechter1985generalized} and dominance criteria in labels \citep{Desaulniers2005}  so as to limit the number of nodes expanded at each iteration of $A*$, thus decreasing the computation time of each iteration of DSSR.  This is done by using the cost of the lowest reduced cost path from any given $u,d$ (computed when $M_v$ sets are empty for all customers $v$) to the sink to describe an admissible heuristic \citep{dechter1985generalized} for each node $u,M_1\subseteq M_u,d$.  Using the Bellman-Ford algorithm we compute the shortest path from each node to the sink when all $M_u$ sets are empty.

\subsection{Outline of Paper}
We organize this document as follows. In Section \ref{Sec_lit_rev_route} we review related literature on pricing in column generation for vehicle routing problems. In Section \ref{Sec_VRP_Relax} we discuss CVRP and its solution via Column Generation.  
In Section \ref{Sec_derivLARoute} we derive pricing from an integer linear programming (ILP) formulation for LA routes. This pricing procedure involved a simple shortest path calculation which we then integrate into Decremental State Space Relaxation \citep{righini2008new} to produce elementary routes.  
In Section \ref{Sec_exper} we provide experimental validation of our approach. In Section \ref{Sec_conc} we conclude and discuss extensions to our research. 

We apply our LA route solution to CVRP, but our method is applicable to any combinatorial optimization problem solved by Column Generation (CG) where pricing is an elementary resource constrained shortest path problem.

\section{Literature Review}
\label{Sec_lit_rev_route}

\subsection{Decremental State Space Relaxation/ng-routes}
Decremental State Space Relaxation (DSSR)\citep{righini2008new,righini2009decremental} is an iterative technique used to solve elementary resource constrained shortest path problems. These problems correspond to the pricing problems generated over the course of Column Generation (CG) for vehicle routing problems (and other problems in operations research).  
DSSR alternates between \textbf{(1)} generating the lowest reduced cost path partially relaxing elementarity (and enforcing resource feasibility, meaning that the path is feasible with respect to any resource constraints) and \textbf{(2)} augmenting the constraints enforced so as to prevent the current non-elementary solution from being regenerated. Step \textbf{(1)} produces the path with the lowest reduced cost from a super set of the set of elementary routes; this set decreases in size as DSSR proceeds. Termination of DSSR is achieved when the generated path is elementary at which point this path is guaranteed to be the lowest reduced cost elementary route.  In practice DSSR does not need to generate all such constraints. DSSR encodes constraints by associating each customer $u$ with a set of the other customers called its neighborhood which is denoted as $M_u$. The path generated at a given iteration of DSSR does not include any cycle satisfying the following property: The cycle starts and ends at $u$, and $u \in M_v $ for all intermediate customers $v$.  
Generating such a path in step \textbf{(1)} is tackled as a dynamic programming problem, which can be alternatively solved using labeling algorithms \citep{Desaulniers2005}.  Given a non-elementary path generated in \textbf{(1)}, a cycle is identified; then in step \textbf{(2)} the neighborhood sets of all intermediate customers are augmented to include the starting/ending customer of the cycle. The solution time of the labeling algorithm can grow exponentially as a function of the maximum size of any neighborhood.  Specifically for the Capacitated Vehicle Routing Problem (CVRP) the time complexity of step \textbf{(1)} scales on the order of $|N|d_0\sum_{u \in N}2^{|M_u|}$ where $N$ is the set of customers and $d_0$ is the capacity of a vehicle.  In order to not needlessly expand the neighborhood sizes and hence the time used by DSSR, DSSR at step \textbf{(2)} can select the cycle intelligently.  For example, DSSR can select the cycle where the augmentation of neighborhoods results in smallest increase in $|N|d_0\sum_{u \in N}2^{|M_u|}$. 

Ng-routes are routes not containing cycles localized in space (cycles localized in space are cycles where all customers are spatially close together). The ng-route relaxation \citep{baldacci2011new, baldacci2012recent} solves the master problem (MP) over a super set of elementary routes known as ng-routes. Solving the MP over ng-routes is easier than solving the MP over elementary routes since pricing over ng-routes is less computationally demanding than pricing over elementary routes. The MP over ng-routes is empirically not much looser than the MP over elementary routes \citep{baldacci2011new}.
The ng-route relaxation can be understood as an adaptation of the ideas of DSSR to prevent the generation of routes with spatially localized cycles in CG. 
The ng-route relaxation can be understood as DSSR where we do not grow the set on neighbors but instead initialize the neighbors of each given customer to be the set of customers it is spatially nearby.  

The popular ng-route relaxation was introduced by Baldacci, Mingozzi, Roberti in 2011 and quickly became widely used for solving vehicle routing problems of various kinds \citep{baldacci2011new, baldacci2012recent}. See for example \citep{baldacci2013exact, bartolini2013improved, contardo2014new, gauvin2014branch, spliet2015discrete, andelmin2017exact, pecin2017improved, pecin2017new, breunig2019electric, duman2021branch}. 
Martinelli et al proposed a method of combining DSSR with ng-routes to examine routes with ng-neighbor set sizes up to 64 and to solve CVRP instances with up to 200 customers \citep{martinelli2014efficient}. 

%
\subsection{General Dual Stabilization}
The number of iterations of CG required to optimally solve the MP can be dramatically decreased by intelligently altering the sequence of dual solutions generated \citep{Pessoa2018Automation,du1999stabilized} over the course of CG. Such approaches, called dual stabilization, can be written as seeking to maximize the Lagrangian bound at each iteration of CG \citep{geoffrion1974lagrangean}.  The Lagrangian bound is a lower bound on the optimal solution objective  to the MP that can be easily generated at each iteration of CG. In CVRP problems the Lagrangian bound is the LP value of the restricted master problem (RMP) plus the reduced cost of the lowest reduced cost column times the number of customers. Observe that when no negative reduced cost columns exist, the Lagrangian bound is simply the LP value of the RMP. The Lagrangian bound is a concave function of the dual variable vector. The current columns in the RMP provide for a good approximation of the Lagrangian bound nearby dual solutions generated thus far but not regarding distant dual solutions.  This motivates the idea of attacking the maximization of the Lagrangian bound in a manner akin to gradient ascent. Specifically we trade off maximizing the objective of the RMP, and having the produced dual solution be close to the dual solution with the greatest Lagrangian bound identified thus far (called the incumbent solution).  

A simple but effective version of this idea is the the box-step method of \citep{marsten1975boxstep}, which maximizes the Lagrangian bound at each iteration of CG s.t. the dual solution does not leave a bounding box around the incumbent solution.  Given the new solution, the lowest reduced cost column is generated and if the associated Lagrangian bound is greater than that of the incumbent then the incumbent is updated. The simple approach of \citep{Pessoa2018Automation} takes the weighted combination of the incumbent solution and the solution to the RMP and performs pricing on that weighted combination.  
Du Merle et al formalized the idea of stabilized CG in their 1999 paper of that name \citep{du1999stabilized}. That paper proposed a 3-piecewise linear penalty function to stabilize CG. Ben Amor and Desrosiers later proposed a 5-piecewise linear penalty function for improved stabilization \citep{amor2006proximal}. Shortly after, Oukil et al used the same framework to attack highly degenerate instances of multiple-depot vehicle scheduling problems \citep{oukil2007stabilized}. Ben Amor et al later proposed a general framework for stabilized CG algorithms in which a stability center is chosen as an estimate of the optimal dual solution \citep{amor2009choice}. Gonzio et al proposed a primal-dual CG method in which the sub-optimal solutions of the RMP are obtained using an interior point solver that is proposed in an earlier paper by the first author \citep{gondzio1995hopdm}. They examine their solution method relative to standard CG and the analytic center cutting plane method proposed by Babonneau et al. \citep{babonneau2006solving, babonneau2007proximal}. They found that while standard CG is efficient for small problem instances, the primal-dual CG method performed better than standard CG on larger problems \citep{Gondzio2013New}.


\subsection{Dual Optimal Inequalities}

Dual optimal inequalities (DOI) \citep{ben2006dual} provide easily computed provable bounds on the space where the optimal dual solution to the MP lies. In this manner, the use of DOI reduces the size of the dual space that CG must search over and hence the number of iterations of CG. Dual constraints corresponding to DOI are typically defined over one or a small number of variables and hence do not significantly increase the solution time of the RMP (though exceptions exist \citep{haghani2020smooth}). DOI are problem instance or problem domain specific. One such example of DOI is in problems such as CVRP or cutting stock where the cost of any column is not increased by removing customers (or rolls in the cutting stock problem \citep{cuttingstock,lubbecke2005selected}) from the column.  Thus equality constraints enforcing that each customer is covered at least once (for the CVRP example) can be replaced by inequality constraints since the cost of a route is not increased by removing a customer from the route. In the dual representation, this replacement corresponds to enforcing that the dual variable corresponding to the constraint that a customer must be covered is non-negative.  

For the cutting stock problem, we can swap a roll of higher length for one of lower length without altering the feasibility of a column (alternatively known as a pattern).  Thus, it can be established that the dual variables associated with rolls, when ordered by non-decreasing roll size must be non-decreasing \citep{ben2006dual}.  In the primal representation, these bounds correspond to swap operations permitting a roll of a given length to be swapped for one of smaller length.  

In \citep{haghani2020smooth}, it is observed that the improvement in the objective corresponding to removing a customer from a column in CVRP (and also the Single Source Capacitated Facility Location Problem (SSCFLP)) can be bounded. Thus primal operations corresponding to removing customers from columns are provided.  In the dual representation, these operations enforce that the reduced cost of a column should not trivially become negative if customers are removed from it.  In the case of SSCFLP, for a given column, this property states that the dual contribution to the reduced cost for a given customer (included in the column) is treated as the maximum of the following two values: the dual variable for that customer, or the distance from the customer to the facility of that column.  

In \citep{haghani2020smooth} it is observed that the dual variables associated with constraints in problems embedded in a metric space should change smoothly over that space.  This is because the dual variable associated with a given customer roughly describes how much larger the objective of the LP is as a consequence of the given customer existing.  Thus, nearby customers should not normally have vastly different dual variables.  Specifically  \citep{haghani2020smooth} shows that in CVRP for any given pair of customers $u,v$ where $u$ has demand no less than $v$, the dual variable of $u$ plus two times the distance from $u$ to $v$ is no less than the dual variable of $v$.  In the primal form, any such pair corresponds to slack variables that provide for swap operations between customers.
This is extended in \citep{yarkony_Detour_DOI} to permit routes to cover customers nearby existing customers on the route at low cost.  
%
%
\subsection {Graph Generation}

Graph Generation (GG) \citep{yarkony2021graph} is an approach to solving the MP, which employs a more computationally intensive RMP than in standard CG at each iteration of CG. GG is implemented with the aim of decreasing the number of iterations of CG required. As in other dual stabilization approaches, pricing is unmodified. GG associates each column $l$ generated during CG with a directed acyclic graph. On this graph, every path from the source to the sink corresponds to a feasible column where the total cost along the edges in the path is the associated cost of the column. Edges are associated with vectors that when summed for a given path from source to sink produce the vector in the constraint matrix corresponding to the path. 
The set of columns corresponding to paths from the source to the sink can be understood as columns related to the column $l$ by a problem instance specific measure. This set is referred to as the family of column $l$ and includes $l$.  
%

The construction of graphs is problem domain/problem instance specific. The use of larger graphs, which have the possibility of containing larger and more diverse sets of columns, has the possibility to decrease the number of iterations of CG. However the use  of larger graphs also leads to increasing computational difficulty of the RMP at each iteration of CG.  
%
For problems where pricing (not solving the RMP) is the computational bottleneck, GG outperforms standard CG because GG requires far fewer iterations than standard CG.

GG can be understood as a generalization of \citep{de2002lp}, which introduces a compact formulation called the arc based formulation for bin-packing/cutting stock problems that is exactly as tight as the standard CG formulation. GG extends the use of such arc based models to cover more general classes of problems by introducing graphs associated with subsets of the columns. 
In the context of CVRP, the graph corresponds to the standard compact formulation of CVRP except that there is an ordering of customers (which varies between graphs) and a column is feasible if it does not violate this order.  Thus if $u$ comes before $v$ in the order associated with column $l$ then no route can be generated in the graph associated with $l$ where $v$ comes before $u$.  

\citep{yarkony2021graph} observed that customers that are in similar physical locations should be in similar positions on the ordered list so that a route permitted by the graph can service these customers without leaving the immediate area containing all of these customers prematurely. Thus \citep{yarkony2021graph} generates the ordering associated with column $l$ by iteratively adding each given customer $u$ behind the existing customer in the route $l$ nearest $u$.   

GG can have a computational bottleneck at solving the RMP when there are numerous graphs or the graphs are large (many nodes and edges). To this end, Principled Graph Management (PGM) \citep{yarkony2022principled} constructs a sufficient set of edges s.t. optimization over those edges provides the optimal solution to the RMP. PGM alternates between solving the RMP over a subset of the edges in the graphs, and adding to the subset any edges associated with the lowest reduced cost path from the source to the sink (which is the lowest reduced cost column associated with the graph) in each graph, which can be computed as a simple dynamic program. Alternatively, edges on paths with near lowest reduced cost can be added and such edges can be efficiently identified using shortest path algorithms. PGM terminates when the optimal solution to the GG RMP is produced, which is achieved when no negative reduced cost columns exist in the families associated with the graphs. 

\section{Column Generation for Capacitated Vehicle Routing }
\label{Sec_VRP_Relax}

We now consider the Capacitated Vehicle Routing Problem (CVRP) which is defined informally as follows. We are given a set of customers with integer demands, a number of homogeneous vehicles with common capacity, and a depot embedded in a metric space. We seek to cover the customers with a set of ordered lists (of customers) called routes, which are each serviced by a unique vehicle so as to minimize the total distance traveled while ensuring that no vehicle services more demand than it has capacity. Vehicles start and end at the depot.  

We now describe CVRP formally.  We use $N$ to denote the set of customers, which we index by $u$. We use $N^+$ to denote $N$ augmented with the starting/ending depot denoted $(-1,-2)$ respectively. The starting/ending depot are in the same place physically but are treated with different notation for convenience. We denote the set of all feasible routes with $\Omega$, which we index by $l$, and is typically too large to enumerate much less consider in optimization.  We use $d_0$ to denote the capacity of a vehicle and $D$ to denote the set of integers ranging (inclusive) from $0$ to $d_0$. Each customer is associated with a positive integer demand (the demand of the starting/ending depot is zero).  We describe the mapping of customers to routes using  $a_{ul} \in \{0,1\}$ where $a_{ul}=1$ if route $l$ services customer $u$ for any $u \in N$. 
For any $u\in N^+,v\in N^+,d \in D$ combination we set $a_{uvdl}=1$ for route $l$ IFF the vehicle leaves $u$ with $d$ units of capacity remaining and departs immediately for $v$. Note that $a_{uvdl}$ is only defined for $d_0-d_u\geq d\geq d_v$ since we must have $d_v$ units of capacity remaining to service $v$ and must have used $d_u$ units of capacity to service $u$. Similarly the vehicle can not leave the ending depot or enter the starting depot and hence $a_{uvdl}$ terms are not defined for those cases.  The demand of the starting/ending depot is zero.  The set of $u,v,d$ that exist are denoted $Q$.  For any pair of $u,v$ which lie in $N^+$ the cost to travel from $u$ to $v$ is denoted $c_{uv}$ and is the distance between $u,v$. The cost of a route is denoted $c_l$, which is defined as the total distance traveled. We express $c_l$ formally as follows in terms of $c_{uv}$ and $a_{uvdl}$.
\begin{align}
\label{route_cost}
    c_l=\sum_{\substack{uvd \in Q}}c_{uv}a_{uvdl}\quad  \forall l \in \Omega
\end{align}
We write the necessary/sufficient characterization of feasible routes as follows.   
\begin{itemize}
\item The route  starts at the starting depot.
\begin{align}
\label{start_feas}
    \sum_{\substack{u \in N}}a_{-1,u,d_0,l}=1 \quad \forall l \in \Omega
\end{align}
\item The route ends at the ending depot.
\begin{align}
\label{endCon}
    \sum_{\substack{u \in N \\ d_0 - d_u \geq d \geq 0}} a_{u,-2,d,l}=1 \quad \forall l \in \Omega
\end{align}
\item The route services no customer more than once meaning that $a_{ul} \in \{0,1\}$ where $a_{ul}$ is defined in terms of $a_{uvdl} \in \{0,1\} $ terms as follows.
\begin{align}
    a_{ul}=\sum_{\substack{v \in N^+\\ d_0-d_u\geq d \geq d_v}}a_{uvdl} \quad \forall l \in \Omega, u \in N 
\end{align}
\item A vehicle must leave a customer $u$ with the amount of resource it entered with minus $d_u$.
\begin{align}
\label{last_feas}
    \sum_{\substack{v \in N^+\\(v,u,d) \in Q}}a_{vudl}=
    \sum_{\substack{v \in N^+\\ (u,v,d-d_u) \in Q}}a_{uv,d-d_u,l} \quad \forall l\in \Omega,u \in N,d\geq d_u
\end{align}
As an aside, observe that \eqref{last_feas} makes \eqref{endCon} redundant.  
\end{itemize}
Given our definition of $\Omega$ the standard CVRP LP relaxation is written using decision variables $\theta_l$ where $\theta_l=1$ if route $l$ is used and otherwise $\theta_l=0$.  We write the minimum weight set cover LP relaxation over $\Omega$ given $K$ vehicles as $ \Psi(\Omega)$ below, with dual variables written in [].
\begin{subequations}
\label{primal_master}
\begin{align}
    \Psi(\Omega)=\min_{\theta \geq 0}\sum_{l \in \Omega}c_l \theta_l \label{CVRP_obj}\\
    \sum_{l \in \Omega}a_{ul}\theta_l\geq 1 \quad \forall u \in N \quad [\pi_u] \label{CVRP_cover}\\
    \sum_{l \in \Omega}\theta_l\leq K \quad [-\pi_0] \label{CVRP_packing}
\end{align}
\end{subequations}
In \eqref{CVRP_obj} we minimize the total cost of the routes used.  In \eqref{CVRP_cover} we ensure that each customer is covered (serviced) at least once. In \eqref{CVRP_packing} we ensure that no more than $K$ vehicles are used.  Though an optimal solution exists covering each customer exactly once, the use of an $\geq$ instead of $=$ in \eqref{CVRP_cover} is known to accelerate common Column Generation (CG) solutions to \eqref{primal_master} (which is called the master problem or MP), without loosening the LP relaxation, by acting as a dual optimal inequality \citep{desrosiers2005primer}.  


Since the cardinality of the set of routes $\Omega$ can grow exponentially in the number of customers we cannot trivially solve \eqref{primal_master}. Instead CG \citep{barnprice,Desrochers1992} is employed to solve \eqref{primal_master}. CG constructs a sufficient subset of $\Omega$ denoted $\Omega_R$ s.t. solving \eqref{primal_master} over $\Omega_R$ provides an optimal solution to \eqref{primal_master} over $\Omega$. 
To construct $\Omega_R$, we iterate between \textbf{(1)} solving \eqref{primal_master} over $\Omega_R$, which is referred to as the restricted master problem (RMP) and written as $\Psi(\Omega_R)$ 
and \textbf{(2)} identifying at least one $l \in \Omega$ with negative reduced cost, which are then added to $\Omega_R$. Typically the lowest reduced cost column (member of $\Omega$) is generated.  We write the selection of this column as optimization below using $\bar{c}_l$ to denote the reduced cost of $l \in \Omega$.  
\begin{subequations}
\label{pricing}
\begin{align}
    \min_{l \in \Omega} \bar{c}_l \\
    \bar{c_l}=c_l+\pi_0-\sum_{u \in N}a_{ul}\pi_u  \quad \forall l \in \Omega \label{redCostForm}
\end{align}
\end{subequations}
The operation in \eqref{pricing} is referred to as pricing.  
CG terminates when pricing proves no column with negative reduced cost exists in $\Omega$.  This certifies that CG has produced the optimal solution to \eqref{primal_master}.  
We initialize $\Omega_R$ with columns corresponding to a heuristically generated feasible integer solution or using artificial variables that have prohibitively high cost to use in an optimal solution but can be used to create a feasible solution.  This can be done by creating $|N|$ variables each of which covers a customer $u$ with prohibitively high cost but without using a vehicle.  In Algorithm \ref{basicCG} we describe CG in pseudo-code.

\begin{algorithm}[!b]
 \caption{Basic Column Generation}
\begin{algorithmic}[1] 
\State $\Omega_R\leftarrow $ from user
\label{Zline_rec_input_start}
\Repeat
\label{Zline_outer_start}%
\State  Solve for $\Psi(\Omega_R)$, generating $\theta,\pi$  using \eqref{primal_master}
\label{Zline_solve_FRMP}
\State $l_* \leftarrow \min_{l \in \Omega}\bar{c}_l$
\State $\Omega_R \leftarrow \Omega_R \cup l_*$
 \Until{$\bar{c}_{l_*} \geq 0$}
 \State Return last $\theta$  generated.  \label{ZreturnSol}
\end{algorithmic}
\label{basicCG}
\end{algorithm} 

We now consider pricing as an binary integer linear program where each binary decision variable $x_{u},x_{uvd}$ correspond to variables $a_{ul},a_{uvdl}$ respectively enforcing \eqref{start_feas}-\eqref{last_feas}.  
Below we define $\bar{c}_{uv}$ to be the cost of traveling from $u$ to $v$ minus the additional cost corresponding to dual variables so that  $\bar{c}_l=\sum_{(uvd)\in Q}\bar{c}_{uv}a_{uvdl} ,\quad \forall l \in \Omega $, which facilitates the efficient writing of pricing.  
\begin{subequations}
\begin{align}
    \bar{c}_{uv}=c_{uv}-\pi_v \quad  \forall (uvd) \in Q, v \neq -2\\
    \bar{c}_{uv}=c_{uv}+\pi_0 \quad \forall (uvd) \in Q, v = -2 
\end{align}
\end{subequations}
We write the computation of the lowest reduced cost route below in the form of an integer linear program (ILP), which we annotate after presenting the ILP.  
\begin{subequations}
\label{pricingEqFull}
\begin{align}
    \min_{x \in \{0,1\} }\sum_{\substack{(uvd)\in Q}}\bar{c}_{uv}x_{uvd} \label{ObjPricer}\\
    \sum_{\substack{(-1,v,d) \in Q}}x_{-1vd}= 1 \label{pushCon}\\ 
    \sum_{\substack{(uvd) \in Q}}x_{uvd}\leq  1 \quad \forall u \in N \label{elemnt}\\
    \sum_{\substack{(uvd) \in Q}}x_{uvd}=\sum_{\substack{(v,u,d-d_v) \in Q}}x_{v,u,d-d_v} \quad \forall  v\in N,d_0\geq d\geq d_v \label{flowcon}
\end{align}
\end{subequations}

In \eqref{ObjPricer} we minimize the reduced cost of the generated route.  The constraints in \eqref{pricingEqFull} correspond to \eqref{start_feas}-\eqref{last_feas}.  In \eqref{pushCon} we enforce that the vehicle leaves the starting depot exactly once.  
In \eqref{elemnt} we ensure that a customer is visited no more than once.  In \eqref{flowcon} we enforce that the vehicle leaves each customer it services with the appropriate amount of demand.  

We should note that the solution of \eqref{pricingEqFull} is not typically solved as an ILP but instead treated as a resource constrained shortest path problem and tackled with a labeling algorithm \citep{costa2019} for the sake of efficiency. The resources that must be kept track of by the labeling algorithm are \textbf{(a)} the set of customers visited thus far and \textbf{(b)} the total amount of capacity used.  

We should note that CG can be applied to supersets of $\Omega$, where the MP over such a superset is written as $\Psi(\Omega^+)$ for some $\Omega \subseteq \Omega^+$.
The use of such  super-sets may weaken the LP relaxation in exchange for greater tractability. The MP over these routes is designed so no optimal integer solution uses an infeasible route. Thus Branch and Price \citep{barnprice} can be used to produce an optimal integer solution to CVRP. Such supersets are often constructed to not include routes with short cycles or cycles localized in space, but do not fully enforce elementarity \citep{baldacci2011new,righini2008new,Desrochers1992}; however, these routes do enforce all other constraints such as capacity and time windows (in domains where time windows are used). For example, the ng-route relaxation \citep{baldacci2011new} is empirically not much looser than the original LP relaxation \eqref{primal_master}; however, the ng-route relaxation is dramatically more efficient to solve compared to solving \eqref{primal_master} because of the ease of pricing.  
\begingroup
\setlength{\tabcolsep}{10pt} 
\renewcommand{\arraystretch}{1.5} 
\begin {table}[H]
\begin{center}
\begin{tabular}{|p{1.2cm}|p{3cm}||p{10cm}|}
 \hline
Name & space & meaning\\
 \hline
   $N$ & set & set of customers\\
   $N^+$ & set & set of customers plus the depot which is counted as $-1$ for starting and $-2$ for ending even though they are the same place\\
$d_0$ & scalar & amount of capacity in a vehicle \\
 $d_u$ &scalar & amount demand at customer $u$ \\
$c_{uv}$ & $u \in N^+,v \in N^+$&  distance from $u$ to $v$.\\
$\pi_u$ & $u \in N$ & dual variable associated with customer $u$\\
$\pi_0$ & scalar & dual variable associated with enforcing an upper bound on the number of vehicles used  \\
$\Omega$& set & set of elementary routes.  \\
$a_{ul}$ & $u \in N,l \in \Omega$& $a_{ul}=1$ if route $l$ covers customer $u$\\
$a_{uvdl}$ & $u \in N^+,v \in N^+,d \in D,l \in \Omega$& $a_{uvdl}=1$ if the route $l$ leaves $u$ with $d$ units of demand remaining and travels immediately to $v$. Here $u,v,d$ must also lie in $Q$ which is the set of feasible possibilities of $u,v,d$\\
 \hline
\end{tabular}
\end{center}
\caption {\bf{CVRP Route Notation}} 
\label{basicTab}
\end{table}
\endgroup
 
 \section{Local Area Routes}
\label{Sec_derivLARoute}

In this section we introduce Local Area routes (LA routes) and pricing over these routes. LA routes are a subset of the celebrated ng-routes \citep{baldacci2011new} and a superset of elementary routes. LA routes have properties which allow for pricing to be conducted efficiently. Solving $\Psi(\Omega^{LA})$, where $\Omega^{LA}$ is the set of all LA routes, provides a fast to compute LP relaxation for Branch and Price solvers for CVRP. Similarly LA routes can be used inside Decremental State Space Relaxation (DSSR) \citep{righini2008new,righini2009decremental},  in order to produce a fast exact solver for CVRP when solving $\Psi(\Omega)$. This in turn can be integrated with Branch and Price.

We organize this section as follows. In Section \ref{def_LANG_Routes} we provide a formal definition of LA routes. In Section \ref{Sec_la_ILP} we describe pricing over LA routes as a shortest path computation.  In \ref{Sec_LA_fast} we describe the fast computation of terms used in Section \ref{Sec_la_ILP}. In Section \ref{sec_decrment} we discuss the use of the DSSR with LA routes to solve $\Psi(\Omega)$. 

\subsection{Definition of LA routes in terms of ng-routes}
\label{def_LANG_Routes}
In this section we define LA routes in terms of ng-routes. We express LA routes using the following terms associated with each customer $u\in N$. Let $N_u$ and $M_u$ refer to two classes of neighboring customers of $u$, which we refer to as the LA neighbors and the ng-neighbors respectively. We define $N_u$ to be a subset of the customers in $N$, that are nearby $u$ with regards to spatial position.  We may define $M_u$ similarly except that we may choose to grow the set when used with DSSR \citep{righini2009decremental}.  By convention, neither $M_u$ or $N_u$ contain $u$. Also, by convention the starting/ending depot has no ng-neighbors or LA-neighbors; furthermore, no customer considers the starting/ending depot to be a ng-neighbor or LA neighbor.

To assist in our description of LA routes, we first describe ng-routes. In the case of ng-routes, $M_u$ is the set of customers nearby $u$ in space. We use $u^p_{k}$ to refer to the $k$'th customer in the path $p \in P$ where $P$ is the set of paths that are resource feasible (meaning that the path does not service more demand than $d_0$; starts at the source and ends at the sink) but these paths do not necessarily have to be elementary. Note that a customer $u$ may be serviced more than once in such a path and hence use more demand than $d_u$.   Recall that a path is considered to be elementary if the path does not visit any customer more than once.  A path $p \in P$ is an ng-route if spatially localized cycles are not present in $p$ (recall that a cycle is a sub-sequence of a path that visits the same customer at the start and end of the sub-sequence). Simply put, in ng-routes, every cycle starting and ending at a customer $u$ must contain at least one customer that does not consider $u$ to be an ng-neighbor. We now express this property with related exposition below.  
%
\begin{align}
\label{ng_deg}
    (u=u^p_{k_1}=u^p_{k_2}) \rightarrow \exists k_3 \quad s.t. \quad k_1<k_3<k_2  \quad  \quad  v=u^p_{k_3} \quad u \notin M_v, \forall\{ u \in N,k_1<k_2\}
\end{align}
The premise of \eqref{ng_deg} (left side of the $\rightarrow $ in \eqref{ng_deg} ) is true if the same customer is at indexes $k_1$ and $k_2$ and that customer is $u$. The inference (right side of the $\rightarrow $ in \eqref{ng_deg}) is that there must exist a customer at some index $k_3$ that lies between $k_1$ and $k_2$ and does not consider $u$ to be a neighbor. 

We now describe LA-routes in a manner akin to our description of ng-routes. 
To assist in our description let us define for any path $p \in P$ a set of indices that we refer to as special; we denote this set of indices as $Q^p$ for path $p$. The first such special index on path $p$ is $1$ meaning that it refers to the first customer visited in the route. We use $q^p_{j}$ to denote the index of the $j$'th special index (meaning $q^p_1=1$).
Thus the $j$'th special index is associated with the $q^p_{j}$'th customer visited in the route. Given that $q^p_{j-1}$ is associated with customer $v$, we define $q_{j}^p$ to be the first index after $q^p_{j-1}$ s.t. the associated customer $u$ does not lie in $N_v$.  We write $q_{j}^p$ mathematically below for any $j>1$.
\begin{align}
\label{qpdef}
    q_{j}^p\leftarrow \min_{\substack{k>q_{j-1}^p\\ u^p_k \notin N_v}} k \quad \forall j>1; \mbox{where $u^p_{q_{j-1}}=v$}
\end{align}
The property beyond resource feasibility that a path must satisfy to be an LA route is as follows: every cycle starting and ending at a customer $u$ must contain at least one customer at a special index that does not consider $u$ to be a ng-neighbor. This is in contrast to ng-routes where that customer does not have to be at a special index.  Thus LA routes are a subset of ng-routes and a super-set of elementary routes.  Thus the CVRP LP relaxation over LA routes is no looser and potentially tighter than the ng-route relaxation (given a fixed set of ng-neighbors for each customer). Using \eqref{qpdef} we now describe the property beyond resource feasibility that an LA route must satisfy formally with exposition below.
\begin{align}
\label{ng_deg2}
    (u=u^p_{k_1}=u^p_{k_2}) \rightarrow \exists k_3 \quad \mbox{s.t. } k_1<k_3<k_2, \quad k_3 \in Q^p,\quad  v=u^p_{k_3}, \quad u \notin M_v \quad \quad \forall \{u\in N,k_2>k_1\}
\end{align}
The premise of \eqref{ng_deg2} (left hand side of the $\rightarrow $ in \eqref{ng_deg2} ) is true if the same customer is at indexes $k_1$ and $k_2$ and that customer is $u$.  The inference (right hand side of the $\rightarrow $ in \eqref{ng_deg2}) is that there must exist a customer at special index $k_3$ that lies between $k_1$ and $k_2$ and does not consider $u$ to be a ng-neighbor. 
%
%

Observe that both LA routes and ng-routes may visit the same customer more than once and hence $a_{ul}$ is a non-negative integer but is not necessarily of binary value (0 or 1).  However $a_{uvdl}$ is still binary for both ng-routes and LA routes and the cost and reduced cost are defined with \eqref{route_cost} and \eqref{redCostForm} respectively. 

We now consider an example of a ng-route that is not a LA route. In our example, the set $N_u$ is identical to the set $M_u$ for each $u\in N$, $|N|=12$, and the locations of the twelve customers correspond to the set of positions on a classic analog clock. As a result, both the LA neighbors and the ng-neighbors of $u_k$ are [$u_{k-2}$,$u_{k-1}$,$u_{k+1}$,$u_{k+2}$] applied with modulus 12. Thus the neighbors (ng and LA) of $u_4$ (4 o'clock) are $N_{u_4}=[u_2,u_3,u_5,u_6]$ and the neighbors of $u_1$ are $N_{u_1}= [u_{11},u_{12},u_2,u_3]$. Respecting the properties set for our example, the following route is considered to be a feasible ng-route but is not a feasible LA route: -1,$u_3$,$u_1$,$u_5$,$u_1$,-2. Observe that the set of special indices in our route consists of the index $1$, which corresponds to the customer $u_3$.  

\subsection{LA Routes:  Integer Linear Programming Formulation For Pricing}
\label{Sec_la_ILP}

In this section we describe the computation of the lowest reduced cost LA route ($\min_{l\in \Omega^{LA}}\bar{c_l}$), which we show is a shortest path problem, that can be trivially solved with Dijkstra's shortest-path algorithm.  To achieve this we write pricing as a sequence of integer linear programs (ILPs); the last ILP corresponds to the generation of the lowest reduced cost LA route as a shortest path problem.  To assist in our discussion we define $P_{u,v,M_1,M_2,d}$ to be the set of all elementary paths satisfying the following.
\begin{itemize}
    \item The path starts at $u \in N^+$ and ends at $v \in N^+$, where $u$ and $v$ are not equal and $v$ is not in the set ($N_u \cup M_1)$.
    \item No customers in $M_1$ (which is a subset of $M_u$) are visited.  
    \item All customers in $M_2-M_1$ are visited where $M_2 \subseteq M_v$.  
    \item No customers in $M_v-M_2$ are visited.
    \item The total demand serviced is  $d$ (excluding $v$ and including $u$, where the starting/ending depot has no demand) OR if $v$ is the end depot then the total demand serviced (including $u$) does not exceed $d$. 
    \item All customers visited lie in $N_u$ excluding $u,v$.
\end{itemize}
We use $N_p$ to denote the set of customers in path $p$ excluding $u$ and $v$. For purposes of clarity of communication we use $Z$ to denote the set of $[u,v,M_1,M_2,d]$ terms for which $P_{u,v,M_1,M_2,d}$ is not empty. We index the set $Z$ with $z$. Given any dual solution $\pi$ we define the reduced cost associated with a path $p$ as $\bar{c}_{p}$, which is described below.  

\begin{align}
\label{term_help_core1B}
    \bar{c}_{p}=c_{p}+\pi_0[u=-1]-\sum_{w \in N}[w \in (N_p\cup u) ]\pi_w \quad \forall p \in P_z,z \in Z,z=(u,v,M_1,M_2,d)
\end{align}
%
We now formulate pricing over the set of elementary routes as an ILP. We use the following binary decision variables defined on a directed acyclic graph  with edge set $E$ where there is a source node, sink node, and each unique $u,M_1,d$ combination (where $u\in N,M_1\subseteq M_u,d_0 - \sum_{w \in M_1}d_w \geq d\geq d_u$) is associated with a unique node.
\begin{itemize}
\item We set $x_{i,j}=1$ where $i=-1$ (meaning the source), $j=(u,\{\},d_0)$ to indicate that vehicle starts at the starting depot then travels to customer $u$.  The source can be equivalently written as $(-1,\{ \},d_0)$.

\item We set $x_{i,j}=1$ where $i=(u,M_1,d_1)$, $j=-2$  to indicate that a vehicle arrives at $u$ with $d_1$ units of demand remaining and then visits customers in $N_u-M_1$ with total demand less than or equal to $d_1-d_u$ before returning to the depot.   The sink can be equivalently written as $(-2,\{ \},0)$.  

\item We set $x_{ij}=1$ for $i=(u,M_1,d_1)$, $j=(v,M_2,d_2)$ where $v\notin (N_u\cup u \cup M_1)$ to indicate that the vehicle arrives at customer $u$ with $d_1$ units of demand remaining (just prior to servicing $u$) then travels along some path in $P_{z}$ where $z=(u,v,M_1,M_2,d_1-d_2)$.


\item The $x_{ij}$ terms indicate the $P_z$ sets that produce component paths, which are concatenated to form the generated route.  
The selection of component paths is done using $x_p$ terms, which are defined for each $p \in P_z$,$z\in Z$.  We set $x_{p}=1$ to indicate that we use path $p$ for $p \in P_{z}$.  In cases where for a given $p$ there are multiple $z\in Z$ for which $p\in P_z$ then we create an independent variable $x_p$ associated with each $z$.  We refer to this as replication of edges.  

\end{itemize}

The computation of the lowest reduced cost elementary route is described by the following ILP.
\begin{subequations}
\label{advPricerB}
\begin{align}
    \label{objB}\min_{x \in \{0,1\} }\sum_{z\in Z}\sum_{p \in P_{z}}\bar{c}_px_{p}\\
    \label{laconst1B}\sum_{u \in N}x_{-1,(u,\{ \},d_0)}=1\\
    \label{laconst2B}\sum_{ij \in E}x_{ij}=\sum_{ji \in E}x_{ji} \quad \forall i=(u,M_1,d) \\
    \label{laconst3B}\sum_{\substack{ij \in E\\i=(u,M_1,d_1)\\j=(v,M_2,d_2)\\d_1-d_2=d}}x_{ij}=\sum_{p \in P_{z}}x_p \quad \forall z\in Z \quad z=(u,v,M_1,M_2,d) \\
    \label{laconst4B}\sum_{\substack{z \in Z \\ z = (u,v,M_1,M_2,d)}}\sum_{p \in P_{z}}x_{p}[w \in (N_p\cup u)]\leq 1 \quad \forall w \in N
\end{align}
\end{subequations}
The objective equation in \eqref{objB} seeks to minimize the reduced cost of the generated route. The constraint in \eqref{laconst1B} forces the number of vehicles used by this route to be 1. The constraint in \eqref{laconst2B} ensures that for every selected edge entering a given node (excluding the source and the sink), there is another selected edge leaving that node (so that the $x_{ij}$ terms describe a connected path from source to sink). 
In \eqref{laconst3B} we ensure that the $x_p$ terms, which describe the intermediate customers in localized components of the path in space, are consistent with the $x_{ij}$ terms, which describe the non-localized structure of the path in space.   Observe that each $p \in P_z$ (for any $z \in Z$) is associated with exactly one constraint of the form \eqref{laconst3B} because of the replication of edges.  
The constraint in \eqref{laconst4B} ensures that no customer is serviced more than once. 

 Observe that no optimal solution to \eqref{advPricerB}  uses a member of $P_z$ that has a higher cost that a different member of $P_z$ given that the same customers are serviced. We express this concept using $\Omega_z$. We use $\Omega_{z}$ to denote the set of paths in $P_{z}$ that are of lowest cost given fixed set of customers.  We construct $\Omega_z$ as follows.  Given any $\hat{N} \subseteq N_u$ we select the single lowest cost path in $p\in P_{z}$ for which $\hat{N}=N_p$ (breaking ties arbitrarily). The set of paths in the $\Omega_z$ terms over all $z \in Z$ is alternatively called the set of LA arcs.
The determination of the paths that define $\Omega_{z}$ (over all $z\in Z$) is done exactly once \textbf{(not at each round of pricing)} because the membership of $\Omega_z$ does not change with $\pi$. 
We provide detailed exposition on the efficient determination of the paths in $\Omega_{z}$ (over all $z\in Z$) in Section \ref{Sec_LA_fast}, which does not enumerate all of $P_z$. Using $\Omega_{z}$, we reduce the size of \eqref{advPricerB} by replacing optimization over $P_{z}$ with optimization over $\Omega_{z}$, producing the following equivalent ILP.

\begin{subequations}
\label{advPricerC}
\begin{align}
    \eqref{advPricerB}=\min_{x\in \{0,1\}}\sum_{ z\in Z}\sum_{p \in \Omega_{z}}\bar{c}_px_{p}\\
    \sum_{u\in N}x_{-1,(u,\{\},d_0)}=1\\
    \sum_{ ij \in E}x_{ij}=\sum_{ji \in E}x_{ji} \quad \forall i=(u,M_1,d) \\
    \sum_{\substack{ij \in E\\i=(u,M_1,d_1)\\j=(v,M_2,d_2)\\d_1-d_2=d}}x_{ij}=\sum_{p \in \Omega_{z}}x_p \quad \forall z \in Z; z=(u,v,M_1,M_2,d) \\
    \sum_{\substack{z \in Z \\ z = (u,v,M_1,M_2,d)}}\sum_{p \in \Omega_{z}}x_{p}[w \in (N_p\cup u)]\leq 1 \quad \forall w \in N\label{badCon2}
\end{align}
\end{subequations}


%
Since we seek to generate the lowest reduced cost LA route and not  the lowest reduced cost elementary route, we can ignore \eqref{badCon2} since that constraint enforces elementarity.  Thus, an optimal solution to \eqref{advPricerC} (when ignoring \eqref{badCon2}) would never use a $p \in \Omega_{z}$ other than the minimizer of reduced cost over $\Omega_{z}$. Given any dual solution $\pi$ we use $p^{z}$ to denote the member of $\Omega_{z}$ with lowest reduced cost meaning:
\begin{align}
    p^{z}\leftarrow \mbox{arg}\min_{p \in \Omega_{z}}\bar{c}_p
\end{align}
We now write the selection of the lowest reduced cost LA route as an ILP.  
\begin{subequations}
\label{penUlt}
\begin{align}
    \min_{x\in \{0,1\}}\sum_{z \in Z}\bar{c}_{p^{z}}x_{p^{z}}\\
    \sum_{u\in N}x_{-1,(u,\{\},d_0)}=1\\
    \sum_{ ij \in E}x_{ij}=\sum_{ji \in E}x_{ji} \quad \forall i=(u,M_1,d) \\
    \sum_{\substack{ij \in E\\i=(u,M_1,d_1)\\j=(v,M_2,d_2)\\d_1-d_2=d}}x_{ij}=x_{p^{z}} \quad \forall z \in Z; z=(u,v,M_1,M_2,d) \label{ezFix2}
\end{align}
\end{subequations}
The formulation in \eqref{penUlt} permits us to use \eqref{ezFix2} to replace $x_{p^{z}}$ in the objective.
\begin{subequations}
\label{my_pricer_ez2}
\begin{align}
    \eqref{penUlt}=\min_{x\in \{0,1\}}\sum_{z \in Z}\bar{c}_{p^{z}}(\sum_{\substack{ij \in E\\i=(u,M_1,d_1)\\j=(v,M_2,d_2)\\d_1-d_2=d\\ z=(u,v,M_1,M_2,d)}}x_{ij})\\
    \sum_{u\in N}x_{-1,(u,\{\},d_0)}=1\\
    \sum_{ij \in E}x_{ij}=\sum_{ji \in E}x_{ji} \quad \forall i=(u,M_1,d) 
\end{align}
\end{subequations}
Observe that the ILP described by \eqref{my_pricer_ez2} is identical in form to ILPs describing standard shortest path problems. We can thus solve \eqref{my_pricer_ez2} using the Bellman-Ford algorithm to find the shortest path from -1 and -2 in $E$, since $E$ may have negative weights but no negative weight cycles. Since $E$ is a directed acyclic graph, $E$ has no cycles and therefore no negative weight cycles. If $E$ can be transformed to an equivalent representation with no negative edge weights, we can apply Dijkstra's algorithm instead to find the shortest path. Applying Dijkstra's algorithm to an equivalent representation of $E$ (which produces the same path and the same minimal cost) is desirable due to the improved asymptotic time complexity of Dijkstra's algorithm compared to the Bellman-Ford algorithm. We detail how $E$ can be transformed into an equivalent representation over which Dijkstra's algorithm can be used to produce the shortest path below.

We let $d_{-1},d_{-2}$ be defined to be $0$ and $d_i=d$; where $i$ corresponds to $(u,M_1,d)$. Observe that for any path starting at $-1$ and ending at $-2$ where the edges included are the set $\hat{E}$ the following property holds: $\sum_{ij \in \hat{E}}(d_i+([i=-1]*d_0)-d_j)=d_0$. Let $c_{ij}$ be the shorthand for the weight in front of $x_{ij}$ in \eqref{my_pricer_ez2}. Let us offset $c_{ij}$ where $i=(u,M_1,d_1)$ and $j=(v,M_2,d_2)$ by adding $\eta *(d_i+([i=-1]*d_0)-d_j)$ to $c_{ij}$ where $\eta=-\min_{z \in Z}\frac{\bar{c}_{p^z}}{z_d}$ and $z_d$ is the capacity used in $z$ ($z_d=\sum_{w \in \hat{N}_{p^z}\cup u}d_w$). This addition makes every edge non-negative and the optimal path is identical to the shortest path produced by Bellman-Ford so Dijkstra's algorithm can be used.  The optimal path from source to sink is not modified because the addition of the $\eta$ terms increases the cost of every path from source to sink by exactly $d_0\eta$. By subtracting the term $d_0 \eta$ from the cost of the generated path, the original cost of the path (without the $\eta$ additions) is obtained. Using Dijkstra's algorithm, only the $\Omega_{z}$ terms corresponding to expanded nodes are required to be computed. In comparison, Bellman-Ford algorithm requires the computation of all possible $\Omega_{z}$ terms, making Bellman-Ford far slower computationally.  

Furthermore, when finding the lowest cost path using Dijkstra's algorithm, we can ignore expanding nodes that are considered to be dominated by another expanded node. A node $(u,M_1,d_1)$ is considered to be dominated by another node $(u,M_1,d)$ if the cost to reach $(u,M_1,d)$ from the source is less than the cost to reach $u,M_1,d_1$ from the source and $d_1<d$. Choosing to include $(u,M_1,d_1)$ instead of $(u,M_1,d)$ is costlier and leaves less remaining capacity for the path, meaning that no optimal solution uses $u,M_1,d_1$. We now demonstrate this phenomenon formally. Let $\tau_{(u,M_1,d)}$ refer to the cost of the shortest path (minimal sum of cost terms from start to end) in graph $E$ starting at the source and ending at $(u,M_1,d)$. Note that $\tau$ does not consider the addition of $\eta$ terms. Suppose that the node $(u,M_1,d_1)$ has been previously expanded and Dijkstra's algorithm now seeks to expand a node $(u,M_1,d_2)$ for which $d_2<d_1$ and for which $\tau_{u,M_1,d_1}<\tau_{u,M_1,d_2}$. Clearly, the lowest reduced cost LA-route can not include the node $(u,M_1,d_2)$ 
and hence we need not expand the node.  More advanced dominance criteria can be developed but we used this one in experiments for its simplicity.  
%
%
%
%
\subsection{Fast Computation of Lowest Cost Component Paths}
\label{Sec_LA_fast}
In this subsection we describe the pre-computation of lowest cost component path terms, a process that ensures that pricing using LA routes is computationally fast. This pre-computation is done once prior to the first iteration of CG and never needs to be repeated. This section shows that $c_{uv\hat{N}}$ terms are easy to compute when the size of $N_u$ sets are small.  To assist in pre-computing $c_{uv\hat{N}}$ terms, we define terms of the form $p^u_{v,w,\hat{N}}$, where $v \in \hat{N}$, $w \in {\hat{N}}$, $\hat{N} \subseteq N_u$, and $u \in N$ to be the lowest cost elementary path visiting all customers in $\hat{N}$, starting at $v$, and ending at $w$. The cost associated with a given $p^u_{v,w,\hat{N}}$ term is given by the term $ c^u_{v,w,\hat{N}}$.  

We now consider a dynamic programming solution to compute $p^u_{v,w,\hat{N}}$ terms efficiently. We now define the base cases for this dynamic programming solution. Clearly when $\hat{N}$ contains only $v_1$, the only associated $ c^u_{v,w,\hat{N}}$ term is $c^u_{v_1,v_1,\{v_1\}}=0$.  Note that this is the only case where the arguments in $c^u_{v,w,\hat{N}}$ have $v=w$ since the path must be elementary. For any $v \in N_u,w\in N_u$ where $v \neq w$, $c^u_{v,w,\{ v,w\}}=c_{vw}$ by definition of $ c^u_{v,w,\hat{N}}$ terms.  

We compute $c^u_{vw\hat{N}}$ and $p^u_{v,w,\hat{N}}$ terms for increasing sizes of $\hat{N}$, first for size 3. 
To compute $c^u_{v, w, \hat{N}}$ we optimize over the choice of $y$, where $y \in \hat{N}-(w \cup v)$ for the path which is described as follows. The path begins at $v$, followed directly afterwards by $y$, ends at $w$, and contains $\hat{N}-( [v,y,w])$ in between $y$ and $w$ in the lowest cost order. This is shown below.
\begin{align}
\label{udayan_beat_julian_eq_1}
    c^u_{v,w,\hat{N}}=\min_{y \in \hat{N}-(w\cup v)}c_{v,y}+c^u_{y,w,\hat{N}-v}
\end{align}

Given that $y$ is the minimizer of \eqref{udayan_beat_julian_eq_1}, the path associated with $c^u_{v,w,\hat{N}}$ denoted, $p^u_{v,w,\hat{N}}$, is the concatenation of $[v] $ and $p^u_{y,w,\hat{N}-v}$. We write the computation of $c^u_{v,w,\hat{N}}$ and $p^u_{v,w,\hat{N}}$ as a dynamic program in Alg \ref{gen_all_dist}.
As an aside, we should note that the values of $c^u_{vw\hat{N}}$ terms do not change for different values of $u$. Thus Alg \ref{gen_all_dist} may choose to not compute $c^u_{vw\hat{N}}$ terms already computed for different values of $u$. 

We use $p^u_{u,w,\hat{N}}$ to denote the lowest cost elementary path starting at $u$, ending at $w$, and visiting all customers in $\hat{N} \cup u$ where $\hat{N}\subseteq N_u$ and $w \in \hat{N} \cup u$ (in order to ensure elementarity, $w$ can only equal $u$ when $\hat{N}$ is empty).  The associated cost for a given $p^u_{u,w,\hat{N}}$ term is  $c^u_{u,w,\hat{N}}$.  Here $c^u_{u,u,\{ \}}=0$ and $c^u_{u,w,\{w\}}=c_{uw}$.  For the remaining cases, we compute $c^u_{u,w,\hat{N}}$ by conditioning on the second customer (which we  write as $v$) in the path as follows.
\begin{align}
\label{minMe2}
    c^u_{uw\hat{N}}=\min_{v \in \hat{N}-w}c_{uv}+c^u_{vw\hat{N}} \quad \forall u \in N, w \in \hat{N},\hat{N}\subseteq N_u, |\hat{N}|\geq 3
\end{align}
The associated path $p^u_{uv\hat{N}}$ equals $[u,p^u_{vw\hat{N}}]$, where $v$ is the minimizer of \eqref{minMe2}.  We are now able to compute the $c_{uv\hat{N}}$ terms as follows for non-empty $\hat{N}$ sets by conditioning on the customer positioned immediately before $v$.  
\begin{align}
\label{minMe3}
    c_{uv\hat{N}}=\min_{w \in \hat{N}}c^u_{uw\hat{N}}+c_{wv} \quad \forall u \in N,v \in N^+-(N_u\cup u),\hat{N}\subseteq N_u, |\hat{N}|\geq 1
\end{align}
The associated path $p_{uv\hat{N}}$ equals $[p^u_{uw\hat{N}},v]$, where $w$ is the minimizer of \eqref{minMe3}. 

To compute all $\Omega_z$ set terms, we simply iterate over all $u,v,\hat{N}$ terms and place $p_{uv\hat{N}}$ in the appropriate $\Omega_z$. Given $u,v,\hat{N}$, we iterate over $M_1\subseteq M_u$ s.t. $M_1 \subseteq (M_u-(\hat{N}\cup v))$; we add $p_{uv\hat{N}}$ to the $\Omega_z$ set where $z=[u,v,M_1,M_2,d]$ for $d=d_u+\sum_{u \in \hat{N}}d_w$ and $M_2=M_v\cap (M_1 \cup \hat{N} \cup u)$.
\begin{algorithm}
\caption{Fast Computation of $c^u_{v,w,\hat{N}}$ terms} \label{gen_all_dist}
\begin{algorithmic}[1]

\For {$f=3:\max_{u \in N}|N_u|$} \label{256}
\For{$u\in N $}
\For{$\hat{N}\subseteq N_{u}$ s.t. $|\hat{N}|=f$}
\For{$v\in \hat{N}$} 
\For{$w \in \hat{N}-v$}
\State $y\leftarrow \mbox{arg}\min_{y \in \hat{N} -(v\cup w)}c_{vy}+c^u_{y,w,\hat{N}-v}$
\State $c^u_{v,w,\hat{N}}\leftarrow c_{vy}+c^u_{y,w,\hat{N}-v}$
\State $p^u_{v,w,\hat{N}}\leftarrow [v, p^u_{y,w,\hat{N}-v}]$
\EndFor
\EndFor
\EndFor
\EndFor
\EndFor

\end{algorithmic}
\end{algorithm}

\begingroup
\setlength{\tabcolsep}{10pt} 
\renewcommand{\arraystretch}{1.5} 
\begin {table}[H]
\begin{center}
\begin{tabular}{|p{1.5cm}|p{3cm}||p{7cm}|}
 \hline
Name & Space & Meaning \\
 \hline
 $N_u$ & $N_u \subseteq \{N-u\}$ & The LA neighbors of $u$.  
 \\
 $M_u$ & $M_u \subseteq \{N-u\}$ & The ng-neighbors of $u$.
 \\
 $z$ & $z \in Z$ & $z$ is a member of set $Z$ where $Z$ corresponds to the space of $u \in N^+, v^+ \in N, M_1 \subseteq M_u, M_2 \subseteq M_v, d = d_1 - d_2$. 
 \\
 $i=(u,M_1,d)$& $u \in N, M_1 \subseteq M_u,d_0 - \sum_{w \in M_1}d_w \geq d\geq d_u$ &  Being at $i$ means that the nascent route is currently at $u$, has $d$ units of capacity remaining (prior to servicing $u$), and all customers in $M_1 \subseteq M_u$ have been visited and serviced at least once already. \\
  $\Omega_z$ & $z \in Z$&  The set of minimum cost elementary paths covering all possible sets of customers for elementary paths meeting the constraints set by $z=(u,v,M_1,M_2,d)$. \\
 $p^z$ & $\Omega_z$ & The lowest reduced cost path in $\Omega_z$.\\
 $c_{p^z}$ & $p \in \Omega_z,z \in Z$ & The cost of the lowest reduced cost path in $\Omega_z$.\\
 \hline
\end{tabular}
\end{center}
\caption {\bf{Local Area Routes Notation}} 
\label{tab_LANG}
\end{table}
\endgroup

\subsection{Decremental State Space Relaxation}
\label{sec_decrment}
In this section we consider the use of the DSSR \citep{righini2009decremental} to generate the lowest reduced cost elementary route, by exploiting the efficiencies of LA routes.  
We describe the procedure to use DSSR alongside LA routes as follows. We initialize the neighborhoods $N_u$ for all $u\in N$ to be the composed of the nearest customers to customer $u$.  We initialize $M_u=\{ \}$  $\forall u \in N$, but these $M_u$ sets increase in size over iterations of DSSR. We then iterate between the following two steps until the lowest reduced cost LA route generated is elementary.  
\begin{enumerate}
    \item   Solve for the lowest reduced cost LA route, which is denoted $l$.
    \item 
    Find the shortest (by length) cycle in the selected route $l$. The length of a cycle is defined to be the number of special indices between the two repeated customers at the start and end of the cycle.  Consider that the shortest cycle starts and ends with customer $u$ at indexes $k_1,k_2$ where $k_2>k_1$. Now, for each intermediate customer $u^l_{k_3}$ corresponding to a special index at index $k_3$ for which $k_2>k_3>k_1$, add $u$ to $M_{u^l_k}$.  The number of nodes $u,M_1,d$ in the pricing graph grows exponentially in $|M_u|$.  Hence to accelerate step \textbf{(1)} we seek to select the cycle to use intelligently so as to decrease the number of nodes/edges added to the pricing graph.  Instead of using the shortest cycle to gradually increase the sizes of $M_u$ sets, we can alternatively use the cycle for which the corresponding increase in the number of nodes or edges in pricing graph (in step \textbf{(1)}) is least (in our experiments we select the cycle that minimizes the number of nodes added to the pricing graph). 
    
    Let $N^*$ be the set of $N$ for which $M_u$ was augmented. Given the updated pricing graph the $\Omega_z$ terms that need to be updated are restricted to the subset defined as follows:  $\Omega_z$ for $z=[u,v,M_1,M_2,d]$ s.t. $u\in N^*$ or $v \in N^*$.  Note that this update does not involve recomputing the $p_{uv\hat{N}}$ terms or the corresponding $c_{uv\hat{N}}$ terms.  
\end{enumerate}
We now consider the use of A* \citep{dechter1985generalized} to enhance the speed of convergence of DSSR by limiting the number of nodes expanded during Dijkstra's algorithm. Consider that we have a heuristic providing a lower bound on the distance from each node to the sink (also known as the ending depot); this heuristic is consistent, meaning that accuracy of the heuristic for the children is no less than that of the parent. Then, we can choose to expand the node for which the distance from the source (also known as the starting depot) plus the heuristic is minimized.    
Observe that each iteration of DSSR solves a similar problem since the values of dual variables remain the same and very few nodes and edges are created from larger $M_u$ sets. Given empty $M_u$ sets (as is the case during the first iteration of DSSR), we can easily compute the shortest distance from each node $(u,\{ \},d)$ (for customer $u$ and demand remaining prior to servicing $u$ as $d$) to the sink. We denote this heuristic as $h_{ud}$ and refer to the graph where all ng-neighbor sets are empty as the initial graph.  We associate each $h_{ud}$ term to nodes of the form $u,M_1,d$ for each $M_1 \subseteq M_u$. The $h_{ud}$ terms are computed exactly once for each call to pricing and not for each iteration of DSSR. Observe that the $h_{ud}$ terms provide an admissible heuristic for A* since the set of possible paths (in terms of customers on the path) on the graph where all ng-neighbor sets are empty is a superset of the set of paths on any graph in DSSR (for paths starting at the source or any $u,M_1,d$ and ending at the sink). By using an admissible heuristic we ensure that A* provides the optimal solution to the shortest-path problem. We can easily and jointly compute $h_{ud}$ for each $u,d$ via the Bellman-Ford algorithm.  

In CG we need to only produce a negative reduced cost column at each iteration of pricing to ensure an optimal solution to the master problem. In the early iterations of CG, the dual values may not reflect the final dual values and hence exact pricing may not produce columns used in the final solution. Therefore, excessive time spent on pricing at this stage is of limited value, thus motivating the use of heuristic pricing. The use of DSSR inside LA routes motivates the following mechanism to efficiently produce negative reduced cost columns. We map the non-elementary route generated at each iteration of DSSR to an elementary route. If this route has negative reduced cost we return the associated column to the RMP; otherwise, we continue with DSSR. In order to generate this elementary route, we remove each customer that is included more than once (after its first inclusion). For example, given a route $[-1,u_1,u_2,u_3,u_1,u_5,u_2,u_1,-2]$, we would produce the route $[-1,u_1,u_2,u_3,u_5,-2]$ instead since we would remove the second copy of $u_2$, as well as the second and third copies of $u_1$.

\section{Experimental Validation}
\label{Sec_exper}
In this section we demonstrate the value of LA routes to accelerate the solution of the set cover formulation over elementary routes in the Capacitated Vehicle Routing Problem (CVRP).  We quantify this with and without dual stabilization (via Graph Generation with Principled Graph Management \citep{yarkony2021graph,yarkony2022principled}).  We consider a data set of CVRP instances that vary by numbers of customers,  customer demands, and vehicle capacity.
In these experiments, pricing is solved using Decremental State Space Relaxation \citep{righini2008new,righini2009decremental} (DSSR) over LA routes as described in Section \ref{sec_decrment}.  We compare the total running time as well as the total time spent on pricing

We organize this section as follows.  In Section \ref{exper_sub_alg} we consider the algorithms compared.  In Section \ref{exper_sub_imp} we consider the implementation details.  In Section \ref{exper_sub_data} we describe the our data sets of problem instances.  In Section \ref{exper_sub_res} we provide the timing comparisons. In Section \ref{exper_sub_anal} we analyze the results of our experiments.  

\subsection{Algorithms Compared}
\label{exper_sub_alg}
To determine the value of stabilization we consider unstabilized column generation (UCG) and GG+PGM.  We follow the GG+PGM methodology from  \citep{yarkony2022principled}. We varied the number of nearest neighbors in the LA routes relaxations in order to determine the value of larger LA neighbor sets.  

To determine the value of larger LA neighbor sets, we consider the use of LA neighbor sets of size [0,5,10] associated with each customer.  Note that using $0$ LA neighbors corresponds to using standard DSSR.  Thus, there are a total of six variants under consideration (2 possibilities for dual stabilization times 3 possibilities for the number of LA neighbors).

\subsection{Implementation Details}
\label{exper_sub_imp}
All code is implemented in MATLAB and restricted master problems are solved using the MATLAB linprog solver with default options.  All linear programs are solved from scratch each time.  In future work, we intend to use CPLEX and not solve linear programs from scratch.  All experiments were run on a 2014 Macbook pro running Matlab 2016.

We initialize the columns in the RMP for UCG with one column for each customer describing the cost for servicing that single customer. 
 
\subsection{Problem Instance Data Set}
\label{exper_sub_data}
We first consider the effectiveness of our CG solvers on two large data sets.  

Data Set One:  All customers have unit demand.  The number of customers lies in [20,30,40,60] and the vehicle capacity lies in [4,8,10]. For each problem combination over the set of possible customers and the set of possible vehicle capacities, we have up to ten problem instances. In each case, the locations of the customers/depot are generated uniformly over a two-dimensional grid.  We describe the number of problem instances of each permutation of capacity/number of customers  in the following table.
\begin{center}
\begin{tabular}{||c c c ||} 
 \hline
 Vehicle Capacity &Number of Customers& Number of Problem Instances \\ [0.5ex] 
 \hline\hline
 4   & 20  &  10\\ 
 \hline
4   & 30  &  10\\ 
 \hline
    4  &  40  &  10\\ 
 \hline
    8  &  20  &  10\\ 
 \hline
    8  &  30  &  10\\ 
 \hline
    8  &  40   &  7\\ 
 \hline
    8  &  60   &  1\\ 
 \hline
   10  &  20  &  10\\ 
 \hline
   10  &  30  &   7\\[1ex] 
 \hline
\end{tabular}
\end{center}
Data Set Two:  All customers have integer demand ranging uniformly over the range [1,10]. The vehicle capacity lies in [20,30,40], but the problems are otherwise identical in terms of the random distribution to data set one.
We describe the number of problem instances of each permutation of capacity/number of customers  in the following table.
\begin{center}
\begin{tabular}{||c c c ||} 
 \hline
 Vehicle Capacity &Number of Customers& Number of Problem Instances \\ [0.5ex] 
 \hline\hline
 20  &  20  &  10\\ 
 \hline
   20  &  30  &  10\\ 
 \hline
   20 &   40 &   10\\ 
 \hline
   30  &  20   & 10\\ 
 \hline
   30 &   30  &  10\\ 
 \hline
   30  &  40  &   5\\ 
 \hline
   40  &  20 &   10\\ 
 \hline
   40  &  30 &   10\\[1ex] 
 \hline
\end{tabular}
\end{center}

We also consider the effectiveness of our CG solvers on a smaller data set. We analyze the speed of LP convergence for problems in this data set at the individual instance level.

Smaller Data Set: All customers have integer demand ranging uniformly over the range [1,10]. The capacity of the vehicle is 40. The number of customers lies in [20,30], and only one problem instance is considered per problem combination over the set of possible customers and the set of possible vehicle capacities. 

\subsection{Results}
\label{exper_sub_res}
\begin{figure}[!hbtp]
	\includegraphics[width=0.49\linewidth]{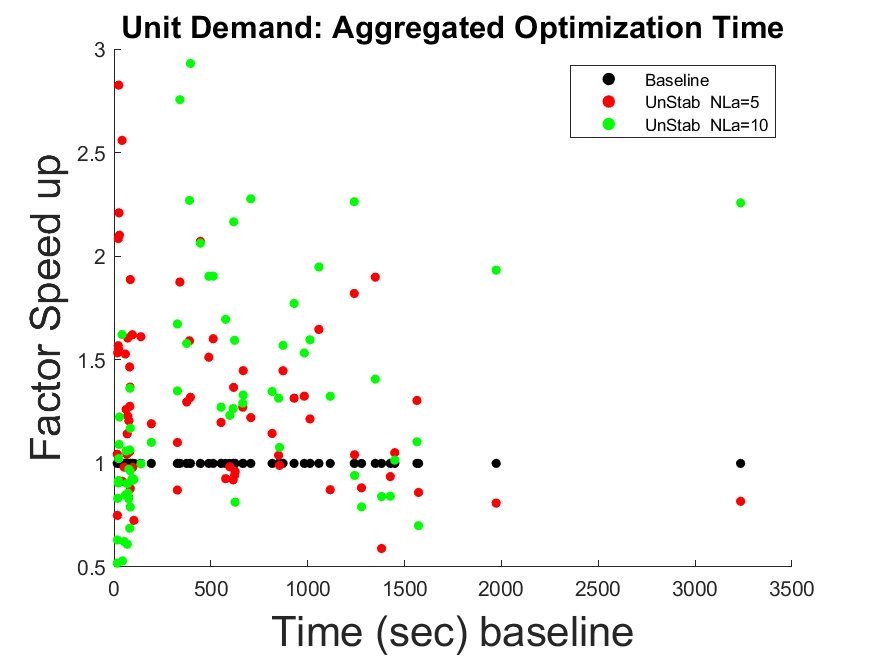}
	\includegraphics[width=0.49\linewidth]{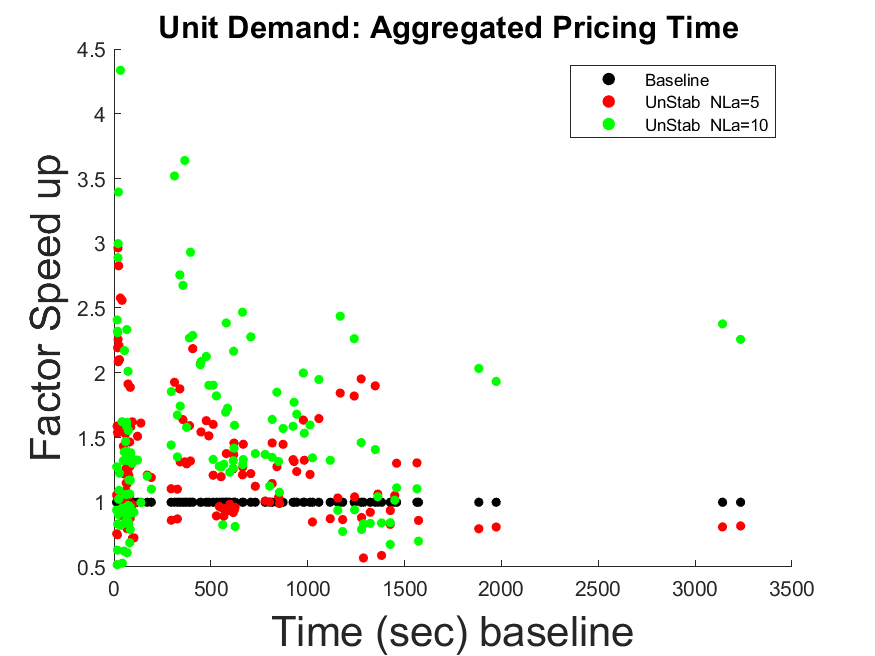}\\
	\includegraphics[width=0.49\linewidth]{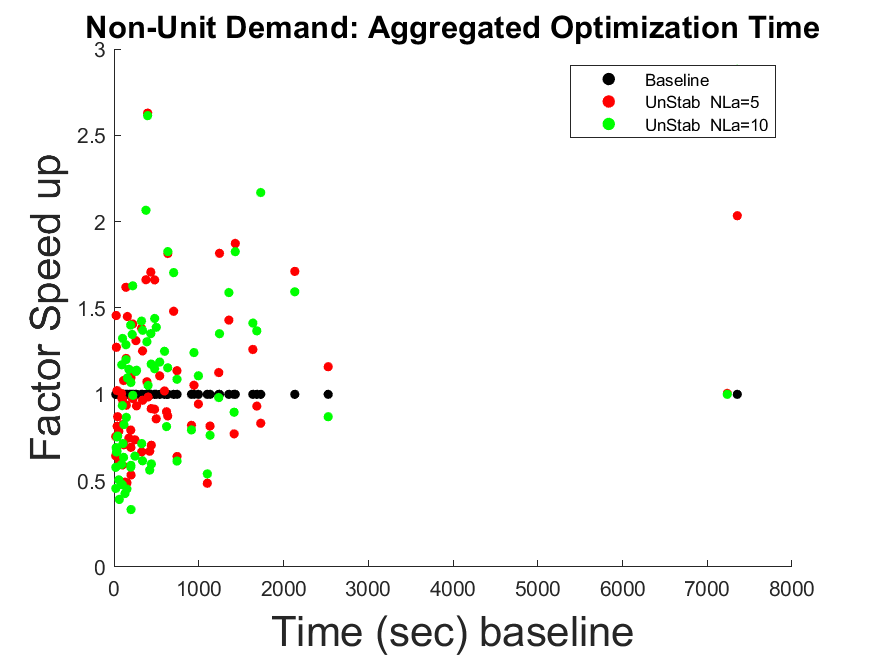}
	\includegraphics[width=0.49\linewidth]{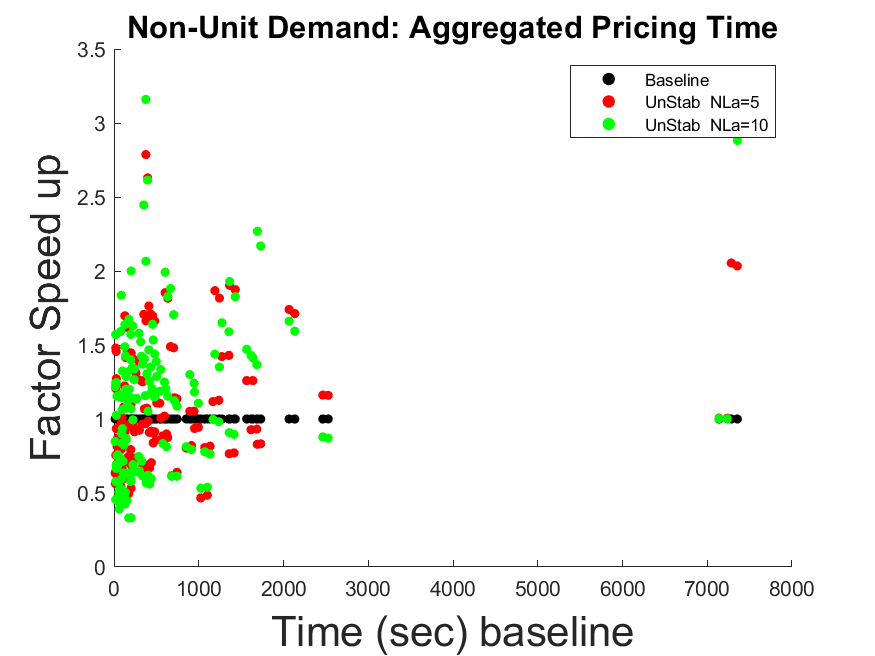}
	\caption{Factor Speed ups for Unstabilized CG solvers (top row unit demand, bottom row, non-unit demand).}
	\label{fig_unit_1}
\end{figure}

\begin{figure}[!hbtp]
	\includegraphics[width=0.49\linewidth]{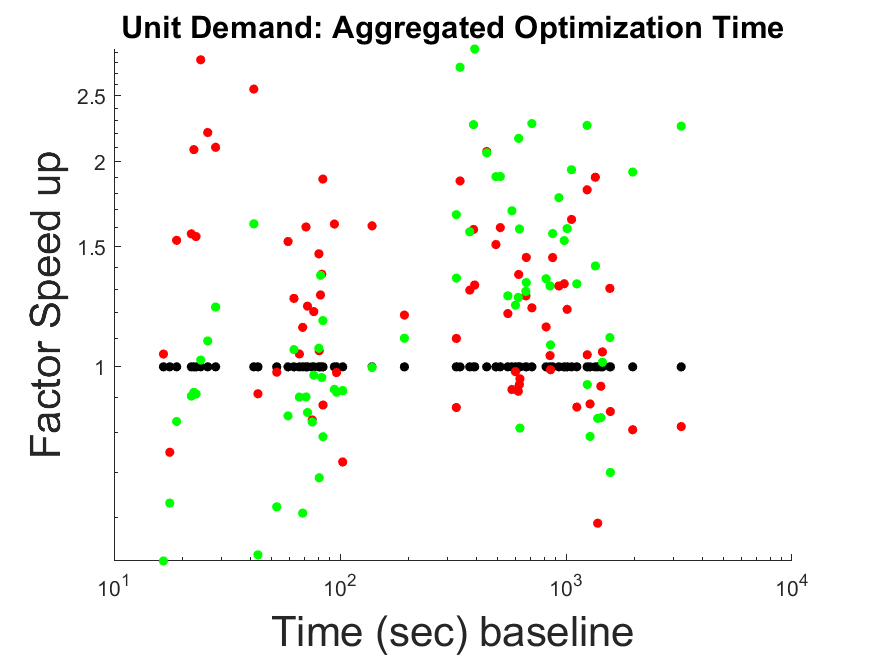}
	\includegraphics[width=0.49\linewidth]{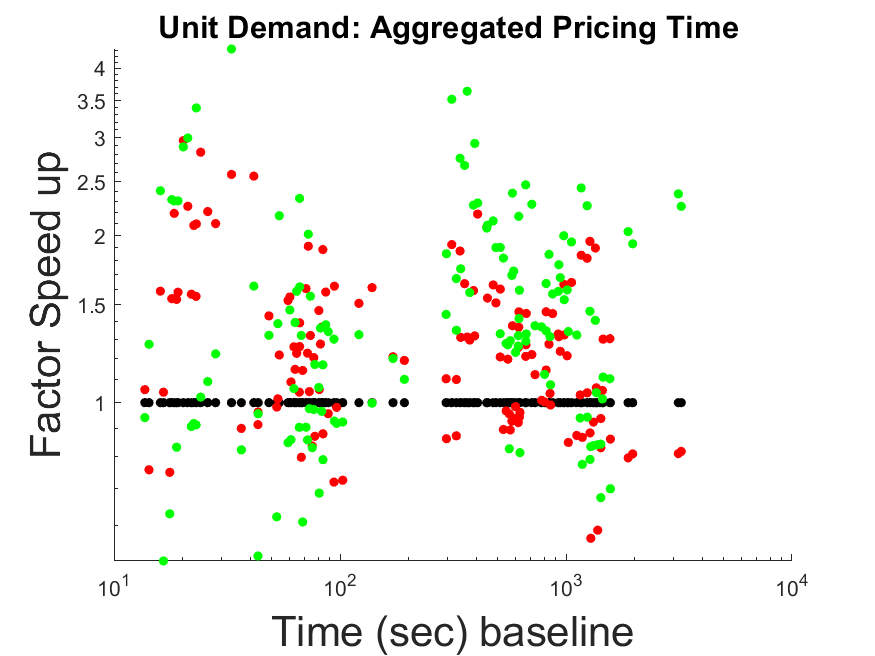}\\
	\includegraphics[width=0.49\linewidth]{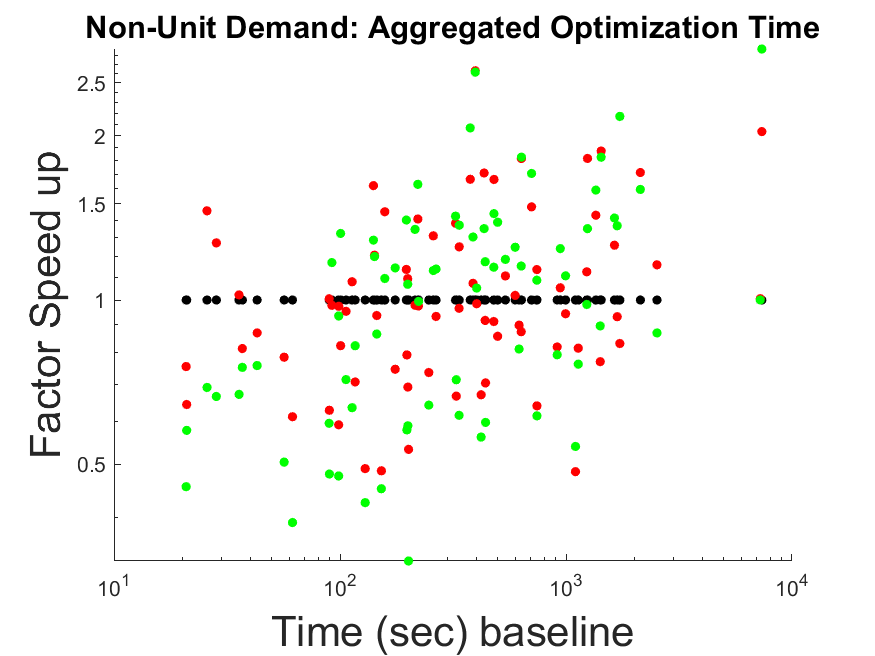}
	\includegraphics[width=0.49\linewidth]{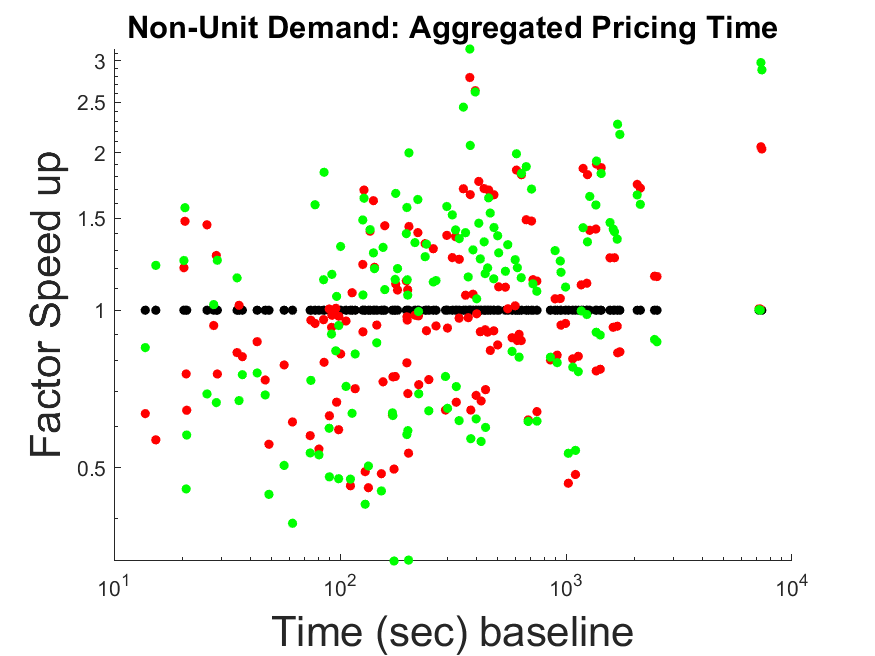}
	\caption{Factor Speed ups for UCG solvers (top row unit demand, bottom row, non-unit demand), but on a log-scale for both the x and y axis. Note for the top row, the y-axis starts at 0.5. The plots in this figure follow the same legend as displayed by plots in Figure \ref{fig_unit_1} (black corresponds to the baseline solver, red corresponds to the UCG solver with LA neighbor set sizes of 5, green corresponds to the UCG solver with LA neighbor set sizes of 10).}
	\label{fig_unit_2}
\end{figure}

\begin{figure}[!hbtp]
	\includegraphics[width=0.49\linewidth]{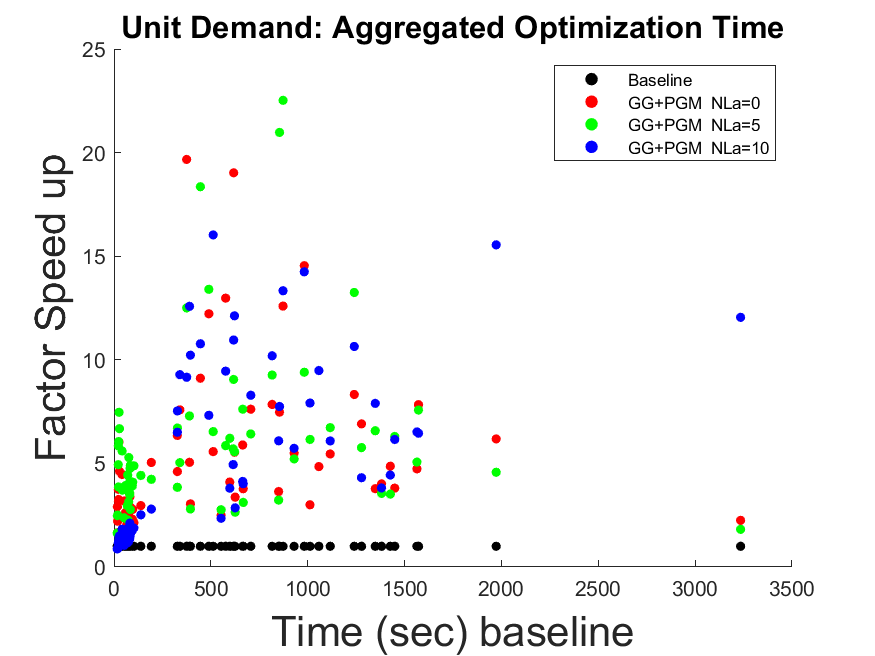}
	\includegraphics[width=0.49\linewidth]{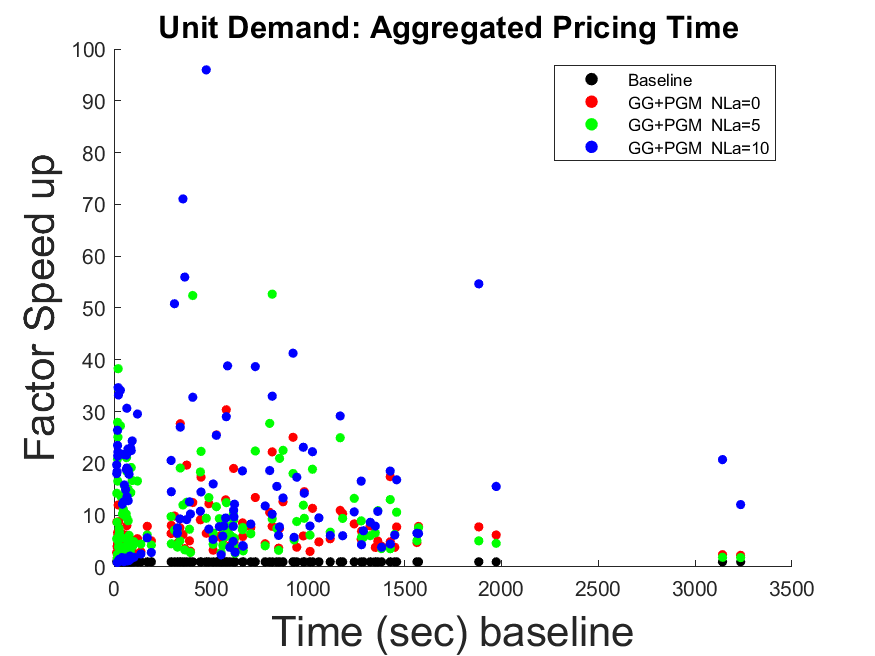}\\
	\includegraphics[width=0.49\linewidth]{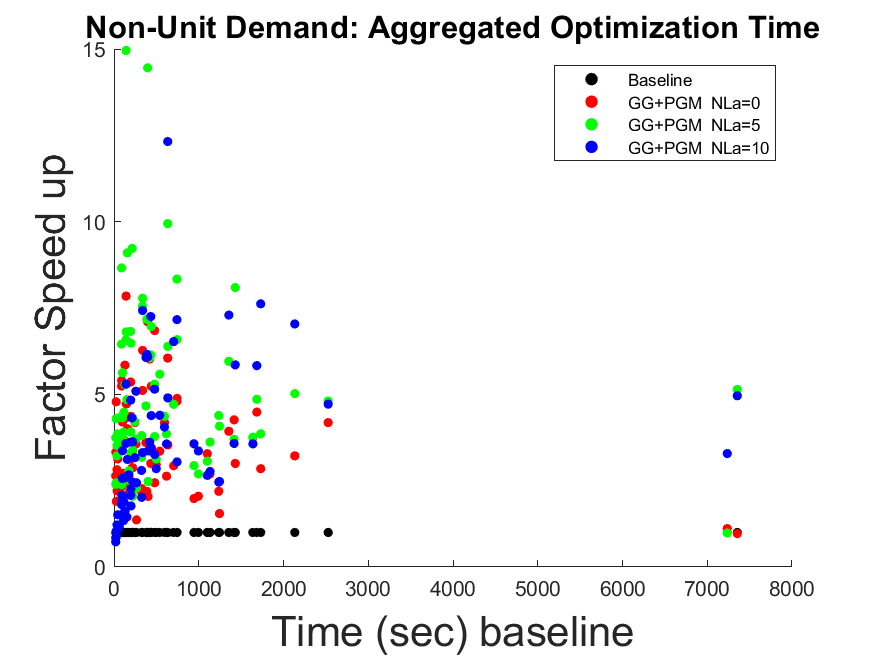}
	\includegraphics[width=0.49\linewidth]{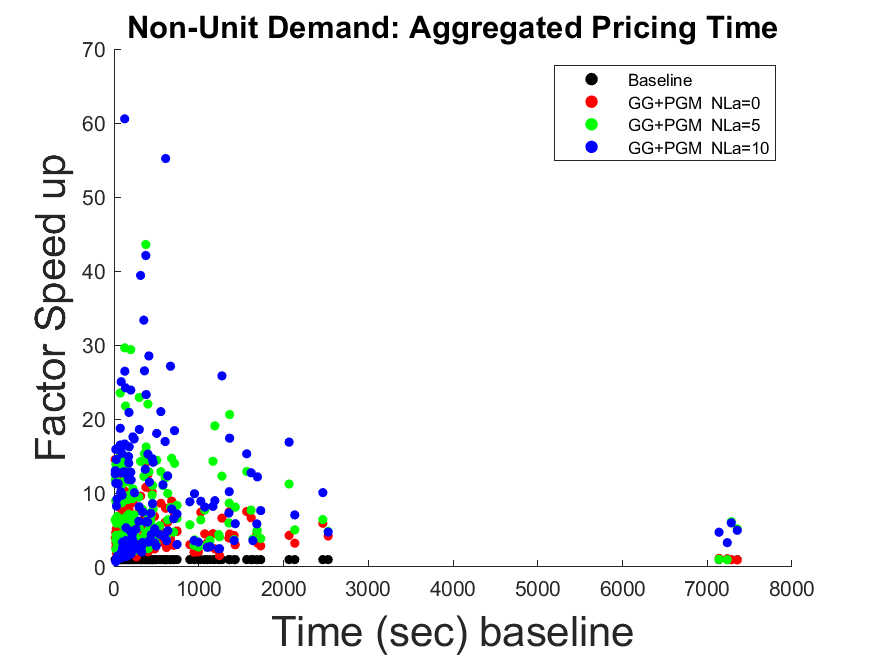}
	\caption{Factor Speed ups for CG solvers incorporating PGM (top row unit demand, bottom row, non-unit demand).}
	\label{fig_unit_3}
\end{figure}

\begin{figure}[!hbtp]
	\includegraphics[width=0.49\linewidth]{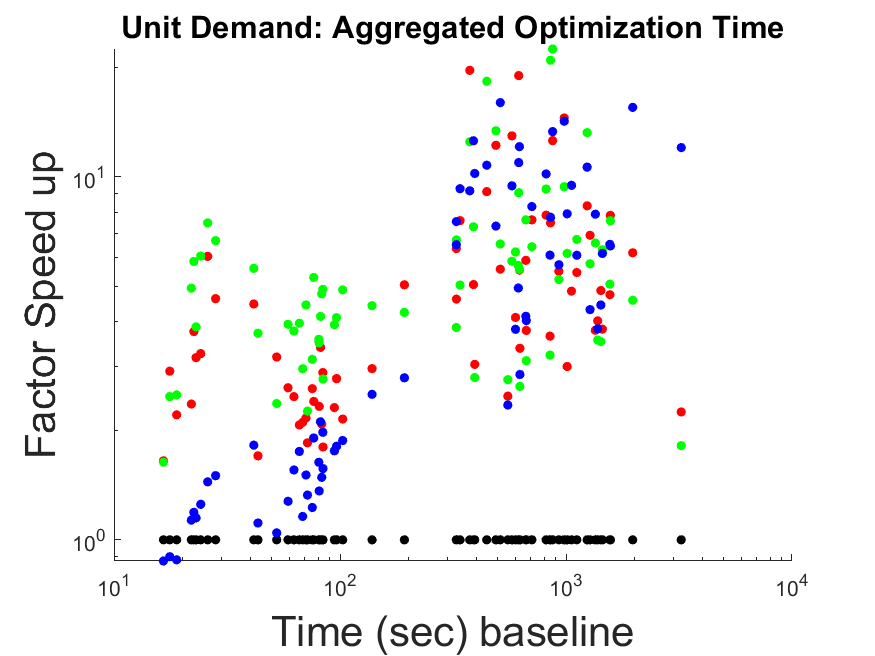}
	\includegraphics[width=0.49\linewidth]{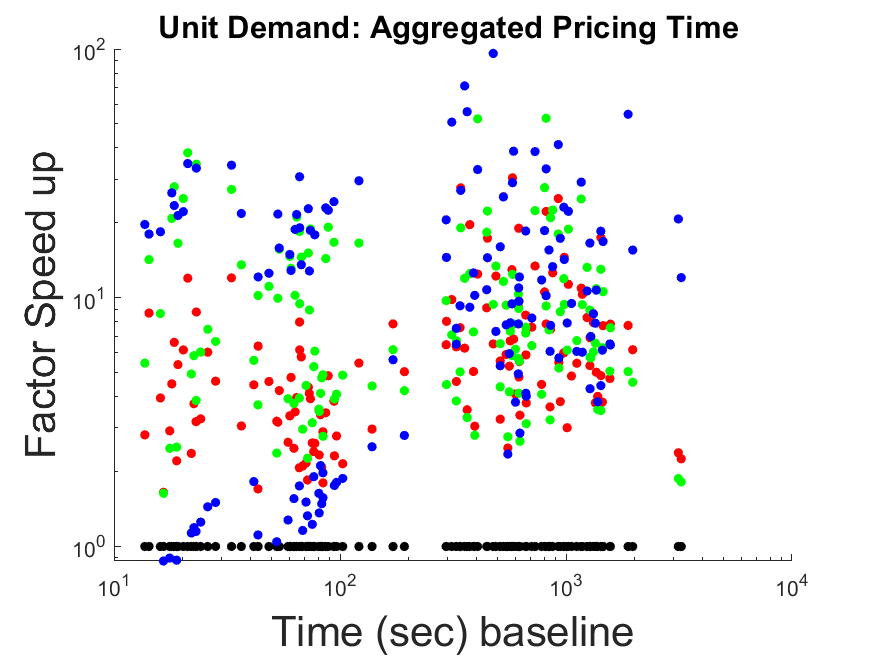}\\
	\includegraphics[width=0.49\linewidth]{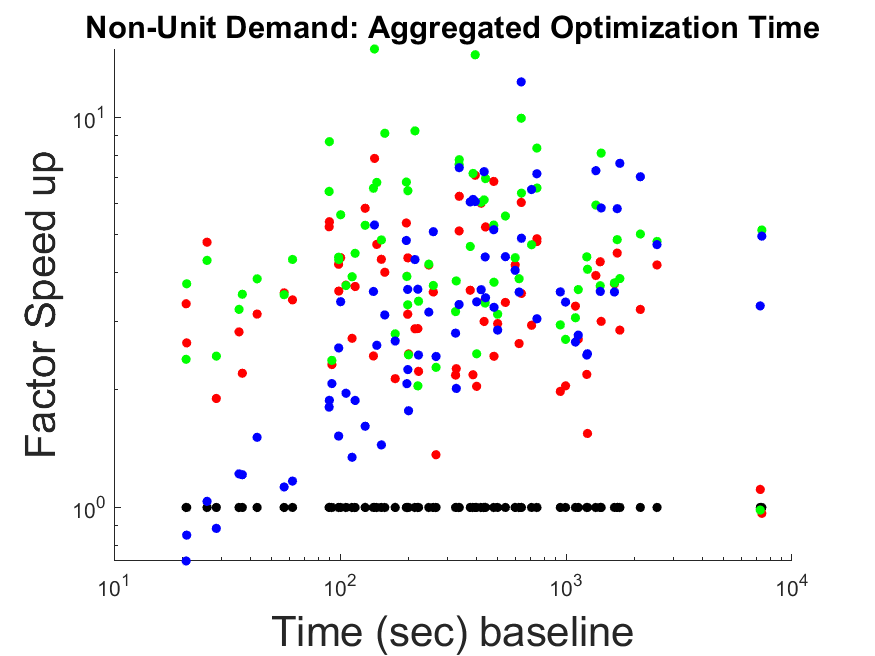}
	\includegraphics[width=0.49\linewidth]{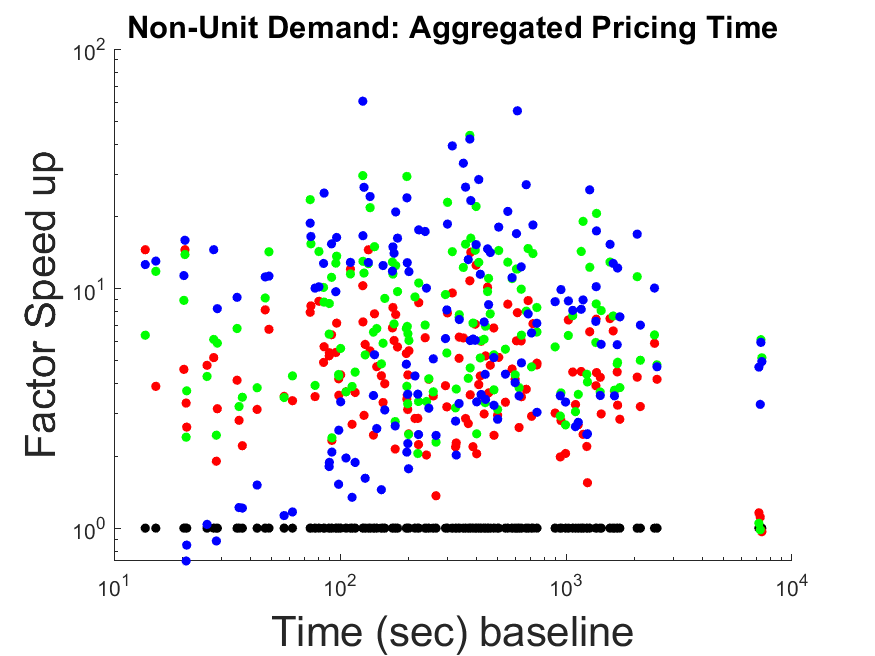}
	\caption{Factor Speed ups for CG solvers incorporating PGM (top row unit demand, bottom row, non-unit demand), but on a log-scale for both the x and y axis. Note for the top row, the y-axis starts at 0.5. The plots in this figure follow the same legend as displayed by plots in Figure \ref{fig_unit_3} (top row unit demand, bottom row, non-unit demand), but on a log-scale for both the x and y axis. Note for the top row, the y-axis starts at 0.5. 
	}
	\label{fig_unit_4}
\end{figure}

For all of our figures, we utilize the  UCG solver with LA neighbor sets of size 0 as our baseline solver which corresponds to standard DSSR. Furthermore, all of our figures represent improvements in computational time taken by CG solvers as a "Factor Speed up" in comparison to the baseline solver. The Factor Speed up of a CG solver for a specific problem is measured by a fraction with the time taken by the CG solver to solve the problem as the denominator  and the time taken by the baseline to solve the problem as the numerator. In Fig \ref{fig_unit_1}, Fig \ref{fig_unit_2}, Fig \ref{fig_unit_3}, and Fig \ref{fig_unit_4}, we have two rows of scatter plots, where the top row shows comparisons in time taken between CG solvers on the first data set, and the bottom row shows the same comparisons on the second data set.  A given data point in all of these figures describes the Factor Speed up of a CG solver (given by the color of the data point) and the time required by the baseline solver for a problem instance; which are encoded on the $x,y$ axis for factor speed up and time required for baseline respectively.
In Fig \ref{fig_unit_1}, we represent the improvements in computational time taken by various implemented  UCG solvers. We represent these improvements on a logarithmic scale for both the x-axis and the y-axis in Fig \ref{fig_unit_2}. 
In Fig \ref{fig_unit_3}, we represent the improvements in computational time taken by various implemented CG solvers utilizing GG+PGM. We represent these improvements on a logarithmic scale for both the x-axis and the y-axis in Fig \ref{fig_unit_4}.

\begin{figure}[!hbtp]
	\includegraphics[width=0.49\linewidth]{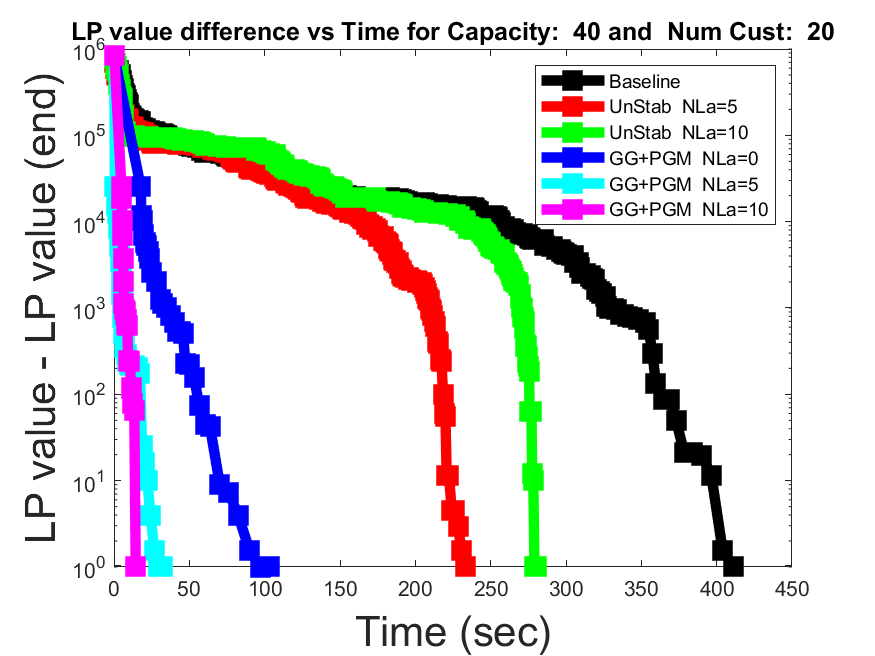}
	\includegraphics[width=0.49\linewidth]{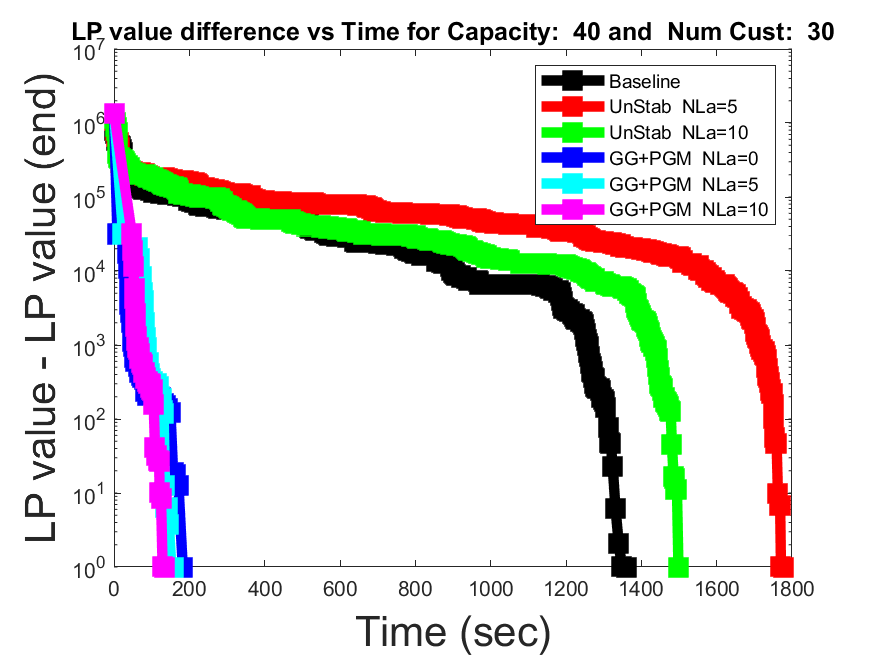}\\
	\includegraphics[width=0.49\linewidth]{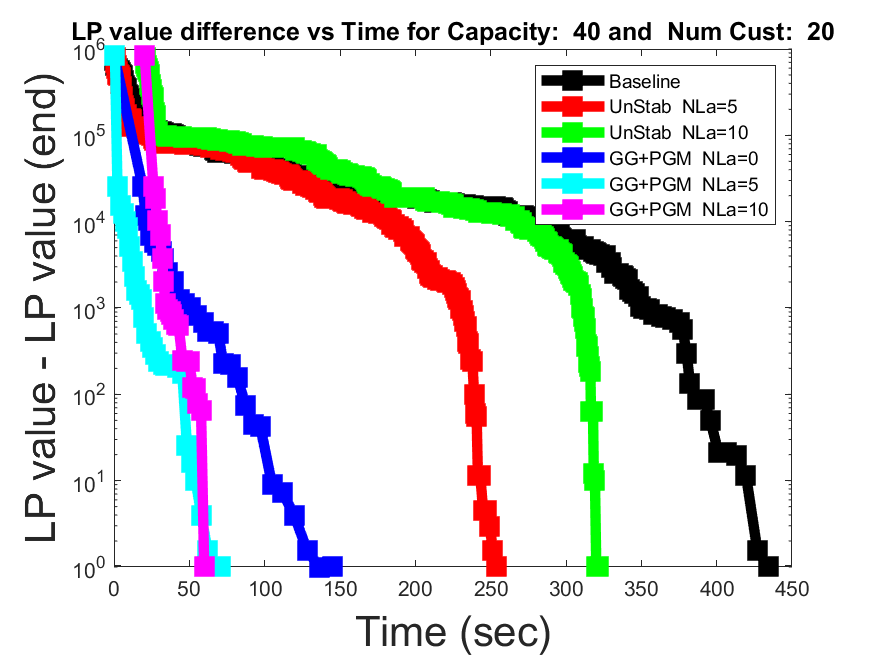}
	\includegraphics[width=0.49\linewidth]{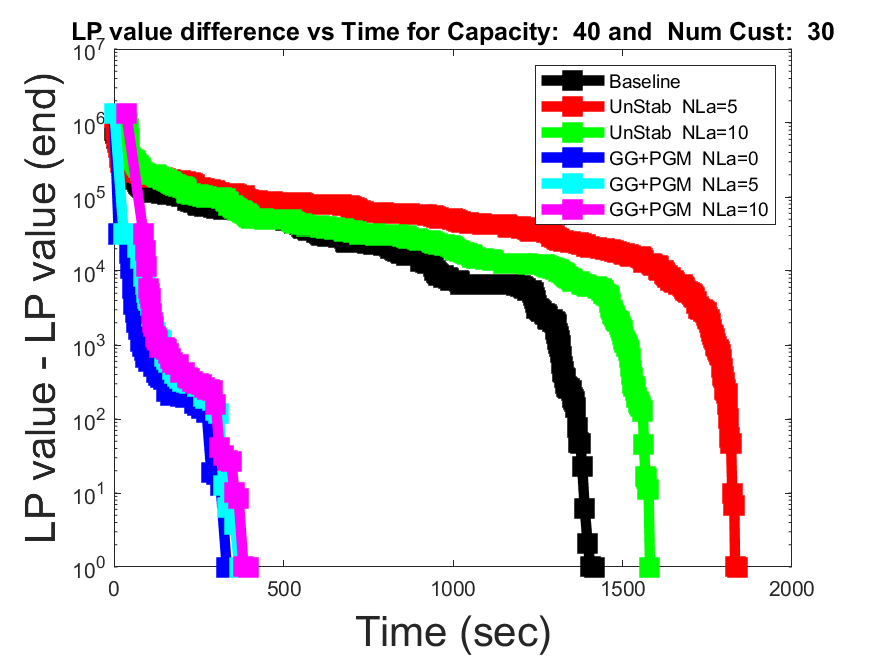}
	\caption{Factors Speed ups for a smaller non unit demand dataset (top row LP convergence rates taking into account time needed for CG, bottom row LP convergence rates only taking into account time needed for pricing, each column corresponds to one problem instance). It is important to note that the metric measuring the difference in LP values always has an addition of 1 to allow for an LP difference of 0 to appear in the graphs.}
	\label{fig_unit_5}
\end{figure}

Data sets 1 and 2 are used alternatively in Tables \ref{speedup_tab_1} and \ref{speedup_tab_2}. These tables describe the proportion of problem instances for which a given approach achieves at least a Factor Speed up for problem instances requiring at 500 seconds to be solved by the baseline solver.  
In Fig \ref{fig_unit_5}, for the first row, we represent the rate of LP convergence for our CG solvers for two problems, taking into account the time taken by CG. In the second row, we represent the rate of LP convergence for our CG solvers for two problems, taking only into account the time taken by pricing. Each column corresponds to LP convergence plots for one of the two problems.
%
In these plots  the x-coordinate shows the computational time taken thus far at a point in time by the CG solver, while the y-coordinate shows the difference in the values of the last given LP solution and the final LP solution on a CVRP problem instance. An addition of 1 is always given to this difference of LP values in order to ensure the graph can always display the y-axis and the data on a log-scale. 
	
\begin{table}

\begin{center}
\begin{tabular}{c|c c c c c c c |} 
 \hline
 Factor Speed up& UCG LA=5& UCG LA=10 & GG+PGM LA=0  & GG+PGM LA=5& GG+PGM LA=10\\
  [0.5ex] 
 \hline\hline
 $*$1   &   0.6543  &  0.8519  &  1.0000  &  1.0000  &  1.0000 \\ 
 \hline
 $*$2&   0.0247 &  0.2346  &  1.0000  &  0.9753  &  0.9877\\ 
  \hline
 $*$5&   0   &    0 &   0.7037  &  0.7531  &  0.8519 \\ 
  \hline
 $*$10& 0    &   0  &  0.2346  &  0.2963  &  0.5185 \\ 
  \hline
 $*$20&   0 &    0  &  0.0617  &  0.0988 &   0.2346\\ 
  \hline
 $*$40&  0  &   0   &   0  &  0.0247  &  0.0741\\ 
  \hline
 $*$60&  0   &    0    &     0    &    0  &  0.0247\\
 \hline
\end{tabular}
\end{center}
\caption{Factor Speed up:  Proportion of problems achieving at least a given speedup over the baseline as a function of the approach used.  We consider problem instances of data set 1 where the baseline solver requires at least 500 seconds to solve the master problem. Only the time taken up during pricing is considered.  A total of 81 problem instances are included. 
\label{speedup_tab_1}
}
\end{table}

\begin{table}

\begin{center}
\begin{tabular}{c|c c c c c c c |} 
 \hline
 Factor Speed up& UCG LA=5& UCG LA=10 & GG+PGM LA=0  & GG+PGM LA=5& GG+PGM LA=10\\
  [0.5ex] 
 \hline\hline
 $*$1   &    0.6250 &   0.7083  &  0.9583  &  0.9792  &  1.0000 \\ 
 \hline
 $*$2&   0.2083  &  0.2917  &  0.9167  &  0.9583  &  1.0000\\ 
  \hline
 $*$5&  0 &   0  &  0.2708  &  0.6250  &  0.6458 \\ 
  \hline
 $*$10& 0   &      0    &     0  & 0.2708 &   0.3333 \\ 
  \hline
 $*$20&   0  &   0     &    0  &  0.0208  &  0.0833\\ 
  \hline
 $*$40&  0   &     0    &   0     &   0  &  0.0208\\ 
  \hline
 $*$60&  0   &     0     &   0   &   0  &   0   \\ 
\\
 \hline
\end{tabular}
\end{center}
\caption{Factor Speed up:  Proportion of problems achieving at least a given speedup over the baseline as a function of the approach used.  We consider problem instances of data set 2 where the baseline solver requires at least 500 seconds to solve the master problem. Only the time taken up during pricing is considered.  A total of 61 problem instances are included. 
\label{speedup_tab_2}
}
\end{table}

\subsection{Analysis}
\label{exper_sub_anal}

We observe large improvements in computation time taken to solve various CVRP problems by solvers combining GG+PGM with LA route relaxations compared to the baseline CG solver. These improvements are generally larger when problem instances which take the baseline longer to solve. We notice these properties in Fig \ref{fig_unit_3}, where the various CG solvers incorporating LA route relaxations with GG+PGM showcase large Factor Speed ups, even for problems which are quickly solved by the baseline CG solver. However, for problems which take larger amounts of time for the baseline CG solver to solve, we often notice larger Factor Speed ups occur more consistently, as shown in the scatter plot measuring the aggregated optimization time for unit demand problem instances. In Fig \ref{fig_unit_4}, this trend is even more apparent, and we can also more clearly see the differences between CG solvers incorporating PGM with different sizes of LA neighbor sets. From Fig \ref{fig_unit_4}, we notice that CG solvers with larger LA neighbor set sizes generally cause greater improvements for larger problems, and the consistency and size of these improvements become more pronounced as the problems increase in size. 

For UCG solvers, we often see sizable improvements in computation time taken for many CVRP problems, but we occasionally see none, or even negative, improvement for some problem instances, as shown in \ref{fig_unit_1}. For larger problems, UCG solvers seem to be more likely to produce larger Factor Speed ups when compared to the baseline solver. Furthermore, in \ref{fig_unit_2}, we notice that the  UCG solver using larger LA neighbor set sizes shows greater Factor Speed ups more consistently and with greater magnitude as problem size increases. Still, in comparison to GG+PGM solvers, the results shown by our UCG solvers are much poorer.

Using Tables \ref{speedup_tab_1} and \ref{speedup_tab_2}, we can more quantitatively gauge both the effectiveness and consistency of our UCG and GG+PGM solvers. In both tables, we notice that the GG+PGM solvers showcase at least Factor Speed ups of 2 in at least 91.67 \%\ of all problem instances, while the UCG solvers showcase at least Factor Speed ups of 2 with far less probability. Furthermore, we see that it is likely to achieve Factor Speed ups as large as 10 using our GG+PGM solvers.


In our results, it is important to note that when only considering time taken by the pricing stage, we notice much greater improvements in computational time taken.  
Therefore, when running the CPLEX LP solver and not solving the RMP LP from scratch at each iteration of CG (or GG+PGM), we expect our results considering only pricing time to more closely resemble our results considering all the time used in the CG process. It is also important to note that these results can vary for different computers and linear programming solver configurations, since hardware specifications and certain settings for the LP solver can result in different properties of the intermediate dual solutions, altering the columns generated during pricing \citep{yarkony2020data}. These properties can result in CG solvers utilizing LA route relaxations to be more or less effective for larger sets of LA neighbors or when GG+PGM is used.


\section{Conclusion}
In this paper we introduce Local Area (LA) route relaxations to improve the tractability and speed of Column Generation (CG) based solvers for large scale set cover/partitioning formulations, where pricing is a elementary resource constrained shortest path problem \citep{costa2019}; this is a common framework found in large-scale transportation problems\citep{barnprice,Desrochers1992}.
LA route relaxations employ LA routes which are routes which do not necessarily need to be elementary but cannot contain spatially localized cycles. LA routes consist of LA arcs, which are elementary resource feasible paths starting at a customer and ending at a customer outside the arc's starting customer's neighborhood (set of customers spatially close to the customer) and visiting intermediate customers that are in the neighborhood of the starting customer.
The only arcs that can be part of an optimal solution to pricing are those which describe the minimum cost path given the starting, ending and intermediate customers. The computation of such arcs is done one prior to any calls to CG as a dynamic program which is tractable when the sizes of neighborhoods is not massive. We used neighborhoods of size 10 efficiently in our experiments.  
We show how LA route relaxations can improve the efficiency of solving one such problem: the Capacitated Vehicle Routing Problem. Our approach is an alternative to the acclaimed ng-route relaxations/decremental state space relaxations (DSSR) \citep{baldacci2011new,righini2009decremental} that permits fast pricing. Pricing stitches together a sequence of LA arcs hence cycles localized in space are easily prevented.
LA route relaxations are no looser and potentially much tighter than ng-route relaxations. 
Furthermore, we show that LA routes can be used alongside DSSR to produce negative reduced cost elementary routes efficiently during pricing.

Additionally in future work, we intend to explore the use of subset row inequalities (SRI) in order to tighten the LP formulation \citep{jepsen2008subset}. In this case, valid inequalities are associated with the LA arcs included in the routes in the restricted master problem (RMP). For example, SRI of size three can be written in a slightly weakened form as follows. We enforce that the number of LA arcs (in routes in the RMP) including two or more members of a set of three customers (excluding the final customer of the arc) cannot exceed one. We can also insert LA arcs into the graphs used in Graph Generation (GG) \citep{yarkony2021graph,yarkony2022principled}) which would permit each GG graph to describe optimal orderings of subsets of customers and express these valid inequalities in the GG RMP without weakening the LP relaxation.   

In future work we intend to explore the use of an alternatively defined master problem that can express all LA routes given the set of LA arcs. We would then solve optimization using column generation, where pricing generates LA arcs. 

We also intend to explore the use of LA routes in the context of time windows. Trivially we could compute for each possible tuple of [start and end time, LA arc] the lowest cost time window feasible path. This path starts and ends at the same customers as the LA arc and visits the associated customers s.t. the path starts at the start time and ends prior to the end time.  Then we could apply LA routes on a graph where nodes describe the tuple of [remaining capacity, current time, ng-neighbors visited, current location]. More efficient mechanisms than enumerating all such time windows for each LA arc would be a fruitful area of exploration.  

\label{Sec_conc}
\bibliographystyle{abbrvnat} 
\singlespacing
\bibliography{col_gen_bib}

\end{document}